\documentclass[preprint,3p]{elsarticle}
\biboptions{sort&compress}

\usepackage[pdftex,colorlinks]{hyperref}
 \usepackage[pagewise,mathlines]{lineno}

\journal{Arxiv}

 \usepackage{tikz}
 \usetikzlibrary{shapes,decorations}

 \usepackage{environ}

\definecolor{mycolor}{RGB}{255,251,204}
\tikzstyle{mybox} = [draw=yellow, very thick,
rectangle, rounded corners, inner sep=5pt, inner ysep=10pt]
\tikzstyle{fancytitle} =[fill=white, text=red]

\NewEnviron{blc}{\begin{tikzpicture}
	\node [mybox] (box){%
		\begin{minipage}{0.95\textwidth}
		\BODY
		\end{minipage}
	};
	\end{tikzpicture}
}

\usepackage{amsfonts}
\usepackage{amssymb}
\usepackage[utf8]{inputenc}
\usepackage{amsmath}
\usepackage{amsthm}
\usepackage{graphicx}%
\usepackage{enumerate}
\usepackage{xcolor}
\usepackage{color,soul}
\usepackage{tikz}
\usetikzlibrary{patterns}
\usepackage{mathtools}

\newtheorem{thm}{Theorem}[section]
\newtheorem{cor}[thm]{Corollary}
\newtheorem{lem}[thm]{Lemma}
\newtheorem{prop}[thm]{Proposition}

\theoremstyle{definition}
\newtheorem{defn}[thm]{Definition}
\theoremstyle{remark}
\newtheorem{rem}[thm]{Remark}

\numberwithin{equation}{section}

\newcommand{\ep}{\varepsilon}
\newcommand{\vep}{\hat{\varphi}_\ep}
\newcommand{\vepx}{\hat{\varphi}_{\ep x}}
\newcommand{\vepxx}{\hat{\varphi}_{\ep xx}}
\newcommand{\yep}{\hat{u}{_\ep}}
\newcommand{\fep}{\hat{h}_\ep}

\newcommand{\R}{\mathbb{R}}				      

\newcommand{\lnum}{ }
\newcommand{\dps}{\displaystyle}	
\newcommand{\intq}{\int_0^T\int_0^1}	
\newcommand{\intw}{\int_0^T\int_{1-\ep}^1}	
	
\newcommand{\into}{\int_0^1}

\newcommand{\n}[2]{\|#1\|_{_{#2}}}

\usepackage[textsize=scriptsize,textwidth=15mm]{todonotes}
\usepackage{marginnote}
\makeatletter
\renewcommand{\@todonotes@drawMarginNoteWithLine}{%
	\begin{tikzpicture}[remember picture, overlay, baseline=-0.75ex]%
	\node [coordinate] (inText) {};%
	\end{tikzpicture}%
	\marginnote[{
		\@todonotes@drawMarginNote%
		\@todonotes@drawLineToLeftMargin%
	}]{
		\@todonotes@drawMarginNote%
		\@todonotes@drawLineToRightMargin%
	}%
}
\makeatother

\allowdisplaybreaks[1] 

\begin{document}

\begin{frontmatter}

\title{\textsc{Boundary null controllability of degenerate heat equation as the limit of internal controllability}}

\author[UFCG]{B. S. V. Ara\'ujo}
\ead{bsergio@mat.ufcg.edu.br}

\author[RCN]{R. Demarque\corref{mycorrespondingauthor}}
\cortext[mycorrespondingauthor]{Corresponding author}
\ead{reginaldo@id.uff.br}

\author[GAN]{L. Viana}
\ead{luizviana@id.uff.br}

\address[UFCG]{Unidade Acadêmica de Matemática, Universidade Federal de Campina Grande, Campina Grande, PB, Brazil}
\address[RCN]{Departamento de Ciências da Natureza,
	Universidade Federal Fluminense,
	Rio das Ostras, RJ, Brazil}
\address[GAN]{Departamento de Análise,
	Universidade Federal Fluminense,
	Niter\'{o}i, RJ, Brazil}

\begin{abstract}
	In this paper, we recover the boundary null controllability for the degenerate heat equation by analyzing the asymptotic behavior of an eligible family of state-control pairs $((u_{\varepsilon}, h_{\varepsilon}))_{\varepsilon >0}$ solving corresponding singularly perturbed internal null controllability problems. As in other situations studied in the literature, our approach relies on Carleman estimates and meticulous weak convergence results. However, for the degenerate parabolic case, some specific trace operator inequalities must be obtained, in order to justify correctly the passage to the limit argument.
\end{abstract}


\end{frontmatter}

\section{Introduction and Statement of the main result}\label{intro}

Take $T>0$, $\alpha \in (0,2)$ and $\varepsilon \in (0,1)$. Let us set
\[
Q:= (0,T)\times (0,1), \ \ \omega_{\varepsilon} :=(1-\varepsilon ,1) \emph{ and } Q_{\varepsilon} :=(0,T)\times \omega_{\varepsilon}.
\]
In this paper, we prove the existence of a family $(u_{\varepsilon}, h_{\varepsilon})_{\varepsilon >0}$, solving 
	\begin{equation}
	\begin{cases}
	u_{\varepsilon t}- \displaystyle \left(x^\alpha u_{\varepsilon x} \right)_x=h_{\varepsilon} \chi_{\omega_{\varepsilon}}, & (t,x)\in Q,\\
	u_{\varepsilon} (t,1)=0,& \text{ in } (0,T),\\
	\begin{cases}
	u_{\varepsilon} (t,0)=0, &\text{if } \alpha \in (0,1),\\
	\text{or}& \\
	(x^\alpha u_{\varepsilon x})(t,0)=0,&  \text{if } \alpha \in [1,2),
	\end{cases} &  t\in (0,T), \\
	u_{\varepsilon} (0,x)=u_0(x),\ u_{\varepsilon} (T,x)=0 & x\in  (0,1),
	\end{cases}
	\label{dcp}
	\end{equation}
with the following property: $(u_{\varepsilon}, h_{\varepsilon})_{\varepsilon >0}$ converges to $(u,g)$, as $\varepsilon \to 0$, where $(u,g)$ solves 
	\begin{equation}
	\begin{cases}
	u_t- \displaystyle \left(x^\alpha u_x \right)_x=0, & (t,x)\in Q,\\
	u(t,1)=g(t),& \text{ in } (0,T),\\
	\begin{cases}
	u(t,0)=0, & \text{if } \alpha \in (0,1),\\
	\text{or}\\
	(x^\alpha u_x)(t,0)=0,& \text{if } \alpha \in [1,2),
	\end{cases},&  t\in (0,T), \\
	u(0,x)=u_0(x),\ u(T,x)=0& x\in  (0,1),
	\end{cases}
	\label{bcp}
	\end{equation}	
	where $u_{0}$ is an initial data taken in a functional space associated to $\alpha$.
	In other words, we will build a suitable family of distributed state-control pairs for the degenerate parabolic equation  \eqref{dcp}, showing that, in some sense, the limiting pair provides the related  boundary null controllability, described in \eqref{bcp}. Above, the systems \eqref{dcp} and \eqref{bcp} are said \textbf{weakly degenerate}, when $\alpha \in (0,1)$, and \textbf{strongly degenerate}, when $\alpha \in [1,2)$ (see \cite{alabau2006carleman} and Remark \ref{boundarycond}).
	
	In \cite{zuazua1988controlabilite}, Zuazua used Lion's Hilbert uniqueness method in order to obtain the internal exact controllability for the wave equation, when the distributed control acts on an appropriate $\varepsilon$-neighborhood of some portion $\Gamma _{0}$ of the boundary. Based on it, in \cite{fabre1992exact}, Fabre proved the exact boundary controllability of the wave equation as the limit of internal controllability, which means that, in the passage to the limit, when $\varepsilon \to 0$, the $\varepsilon$-neighborhood of $\Gamma _{0}$ shrinks to itself. Recently, in \cite{chaves2020boundary}, Chaves-Silva et al. obtained a similar result for the heat equation. In this current work, we are focused on an analogous investigation about the degenerate heat equation case.
	
	Due to the scope we have chosen, we would like to have a brief talk about the  controllability of one-dimensional degenerate problems. To start the discussion, many applied phenomena are closely related to degenerate parabolic equations, calling a notorious attention to their mathematical point of view. In this case, motivated by the properties already known for the uniformly parabolic case, a more complete qualitative literature for degenerate operators is also welcome (see a well-posedness result in \cite{campiti1998degenerate}, for instance). It certainly includes the Control Theory, where much more development is still desired. In one dimension, it seems to us that  \cite{cannarsa2002nulle} and \cite{cannarsa2004persistent} are the two first articles dealing with the controllability of degenerate parabolic equations, which clearly inspired much relevant work since then (see \cite{cannarsa2005null,cannarsa2006null,cannarsa2008carleman, Boutaayamou2018carleman,wang2018carleman,araruna2018stackelberg,fragnelli2018carleman, runmei2019null,mustapha2019algebraic} and the references therein). Up to now, it is undeniable the strength of the Carleman estimate method, because it provides a refined technique that makes the one-dimensional degenerate controllability field well-understood (see  \cite{alabau2006carleman,martinez2006carleman,cannarsa2007null, cannarsa2008controllability,cannarsa2009controllability} and the references aforementioned). This paper intends to contribute in this field of degenerate parabolic PDEs by providing some singular perturbation result like those proved in \cite{fabre1992exact} and \cite{chaves2020boundary}.

Next, we consider some important functional spaces, introduced in \cite{alabau2006carleman}, which are closely related to the initial data $u_0$ in \eqref{dcp} and \eqref{bcp}.

\begin{defn} [Weighted Sobolev spaces]
Consider $\alpha \in (0,1)$, for the \textbf{weakly degenerate case} (WDC), or $\alpha \in [1,2)$, for the \textbf{strongly degenerate case} (SDC).

\begin{itemize}
	\item[(I)] For the (WDC), we set 
	\begin{equation*}
	H_{\alpha}^1:= \{  u\in L^2(0,1);\ u\mbox{ is absolutely continuous in } [0,1],
	x^{\alpha /2}u_x\in L^2(0,1) \mbox{ and } u(1)=u(0)=0\},
	\end{equation*}
	equipped with the natural norm
\[ \|u\|_{H_{\alpha}^1}:=\left( \|u\|_{L^2(0,1)}^2+\| x^{\alpha /2} u_x\|_{L^2(0,1)}^2 \right) ^{1/2} ;\]
\item[(II)] For the (SDC),
\begin{equation*}
H_{\alpha}^1:= \{  u\in L^2(0,1);\ u\mbox{ is absolutely continuous in } (0,1],
x^{\alpha /2} u_x\in L^2(0,1) \mbox{ and } u(1)=0\},
\end{equation*}  
and the norm keeps the same;

\item[(III)] In both situations, the (WDC) and the (SDC), 
\[
H_{\alpha}^2:= \{  u\in H_{\alpha}^1;\ x^{\alpha /2} u_x\in H^1(0,1) \}
\]
with the norm
$\|u\|_{H_{\alpha}^2}:=\left( \|u\|_{H_{\alpha}^1}^2+\|(x^{\alpha /2} u_x)_x\|_{L^2(0,1)}^2 \right) ^{1/2}$. 
\end{itemize}
\end{defn}	

\begin{rem}\label{boundarycond}
    We  use Dirichlet boundary conditions for the (WDC), while the natural boundary conditions for the (SDC) are  Neumann type conditions. The reason is because the notion of trace at the boundary is unavailable for strongly degenerate spaces, as proved in   \cite[Section 17.1 ]{cannarsa2016global}. 
\end{rem}

At this moment, we are supposed to specify which kind of solution for \eqref{dcp} and \eqref{bcp} we are dealing with.

\begin{defn}\label{transp}
Consider $u_0\in L^2(0,1)$, $h_{\varepsilon} \in L^2(Q_{\varepsilon} )$ and $g\in L^2 (0,T)$.

\begin{itemize}
\item[(a)] We say that $u_{\varepsilon} \in L^2(Q)$ is a \textit{\textbf{solution by transposition of }} \eqref{dcp} if, for each $(F,z^T)\in L^2(Q)\times L^2(0,1)$ given, we have
\begin{equation*}
    \intq uF\, dxdt=\intw hz\,dxdt+\into u_0(x)z(0,x)\, dx,
\end{equation*}
where $z$ solves
\begin{equation}\label{adj-transp}
	\begin{cases}
	-z_t- \displaystyle \left(x^\alpha z_x \right)_x=F, & (t,x)\in Q,\\
	z(t,1)=0,& \text{ in } (0,T),\\
	\begin{cases}
	z(t,0)=0, & \text{ (Weak) }\\
	\text{or}\\
	(x^\alpha z_x)(t,0)=0,&  \text{ (Strong) }
	\end{cases},&  t\in (0,T), \\
	z(T,x)=z^T(x), & x\in  (0,1).
	\end{cases}
\end{equation}
\item[(b)] Analogously, we say $u\in L^2(Q)$ is a \textit{\textbf{solution by transposition of}} \eqref{bcp} if, for every $(F,z^T)\in L^2(Q)\times L^2(0,1)$, we have
\begin{equation}\label{sol.bound}
    \intq uF\, dxdt=-\int_0^T g(t)z_x(t,1)\,dt+\into u_0(x)z(0,x)\, dx,
\end{equation}
where $z$ solves \eqref{adj-transp}.
\end{itemize}
\end{defn}
The aim of this paper is to prove the following theorem.

\begin{thm}\label{main}
    Consider $T>0$ and $u_0\in H_\alpha^1$. Then, for each $\varepsilon \in (0,1)$, there exist
    \[
    \displaystyle u_{\varepsilon} \in H^1(0,T;L^2(0,1))\cap L^2(0,T;H_{\alpha}^2)\cap C^0([0,T];H_{\alpha}^1) 
    \emph{  and  } 
    h_{\varepsilon} \in L^{2} (Q_{\varepsilon})
    \]
    such that:
    
    \begin{enumerate}[(a)]
        \item $(u_{\varepsilon} , h_{\varepsilon})$ solves the distributed null controllability problem \eqref{dcp}, in the sense of Definition \ref{transp}(a);
        \item $u_{\varepsilon} \rightharpoonup u$, weakly in $L^2 (Q)$, and $h_{\varepsilon} \rightharpoonup g$, weakly in $L^2 (0,T;H_{\alpha}^{2})$, as $\varepsilon \to 0$. Moreover, 
        \[
        \displaystyle u\in C^0([0,T];(H_\alpha^1)')\cap L^2(Q),
        \]
        \[
        g\in L^2 (0,T)
        \]
        and $(u,g)$ solves the boundary null controllability problem \eqref{bcp}, in the sense of Definition \ref{transp}(b). 
    \end{enumerate}
\end{thm}

\begin{rem}
Theorem \ref{main} means that the boundary null controllability of the degenerate heat equation can be seen as the limit of some proper family of internal controllability problems. However, it seems to us that the same asymptotic behavior is not expected for arbitrary families solving \eqref{dcp}.
\end{rem}


In order to prove Theorem \ref{main}, our approach is based on the strategies of \cite{chaves2020boundary}. In fact, from a Carleman estimate at the boundary, we obtain the optimal observability constant $C=O(\varepsilon ^{-3})$ (see Remark 2.2 of \cite{chaves2020boundary}). This is a meticulous and crucial point, which makes the passage to the limit in \eqref{dcp} possible.
However, for the degenerate parabolic case, some specific trace operator inequalities must be obtained, in order to justify correctly the passage to the limit argument (see Proposition \ref{abalfa} and Corollary \ref{abalfavec}).


The remainder of this paper is organized as follows: in Section \ref{sec-pre}, we present notations and preliminary results, like improved Carleman estimates and some trace operator inequalities. In Section \ref{familia}, we construct a suitable family of state-control pairs $((u_{\varepsilon}, h_{\varepsilon}))_{\varepsilon >0}$ for which we will study the asymptotic behavior. In Section \ref{principal}, we rigorously pass to the limit in \eqref{dcp}, describing all needed weak convergence results. As a consequence, we achieve our main result (Theorem \ref{main}). In Section \ref{aberto}, we present additional comments involving some possible extensions and open questions. In Appendix \ref{app}, we complement the proof of Proposition \ref{car-frnt}, giving technical details omitted in Section \ref{sec-pre}.

\section{Preliminaries}\label{sec-pre}

In this section, we reunite definitions and results which will support the proof of our main result (Theorem \ref{main}). We should emphasize that the whole discussion presented in this paper arises from the well-posedness of \eqref{dcp}, obtained in \cite{alabau2006carleman}, whose statement is given below:

\begin{prop}\label{well-pos}
	Given $u_0\in L^2(0,1)$ and $h\in L^2(Q)$, there exists a unique weak solution $u$ of 
\begin{equation*}
	\begin{cases}
	u_t- \displaystyle \left(x^\alpha u_x \right)_x=h, & (t,x)\in Q,\\
	u(t,1)=0,& t\in  (0,T),\\
	\begin{cases*}
	u(t,0)=0, & \text{ (Weak) }\\
	\text{or}\\
	(x^\alpha u_x)(t,0)=0,&  \text{ (Strong) }
	\end{cases*},&  t\in (0,T), \\
	u(0,x)=u_0(x), & x\in  (0,1),
	\end{cases}
	\end{equation*}
	such that
	$u\in C^0([0,T];L^2(0,1))\cap L^2(0,T;H_\alpha^1)$.
	In addition, if $u_0\in H_{\alpha}^{1}$, then
	\[u\in H^1(0,T;L^2(0,1))\cap L^2(0,T;H_{\alpha}^2)\cap C^0([0,T];H_{\alpha}^1), \]
	and there exists a positive constant $C_T$ such that
	\begin{equation}\label{ineq1}
	\sup_{t\in[0,T]}\left(\n{u(t)}{H_\alpha^1}^2\right)
	+\int_0^T\left(\n{u_t}{L^2(0,1)}^2+\n{(x^\alpha u_x)_x}{L^2(0,1)}^2\right)\\
	\leq C_T\left(\n{u_0}{H_\alpha^1}^2 +\n{h}{L^2(Q)}^2 \right).
	\end{equation}
\end{prop}

We will also state a result of well-posedness for the problem \eqref{bcp}, which the proof is given in \cite{gueye2014exact} for $\alpha\in (0,1)$, but the same argument works for $\alpha\in (0,2)$.
\begin{prop}
	Given $u_0\in (H_\alpha^1)'$ and $g \in L^2(0,T)$, there exists a unique weak solution (by transposition) $u$ of 
\begin{equation*}
	\begin{cases}
	u_t- \displaystyle \left(x^\alpha u_x \right)_x=0, & (t,x)\in Q,\\
	u(t,1)=g(t),& t\in (0,T),\\
	\begin{cases*}
	u(t,0)=0, & \text{ (Weak) }\\
	\text{or}\\
	(x^\alpha u_x)(t,0)=0,&  \text{ (Strong) }
	\end{cases*},&  t\in (0,T), \\
	u(0,x)=u_0(x), & x\in  (0,1),
	\end{cases}
	\end{equation*}
	such that
	$u\in C^0([0,T];(H_\alpha^1)')$. 	In addition, there exists a positive constant $C$ such that
	\begin{equation}\label{ineq-boundary}
	\n{u}{L^\infty(0,T;(H_\alpha^1)')}\leq C\left(\n{u_0}{L^2(0,1)}+\n{g}{L^2(0,T)}\right).
	\end{equation}
\end{prop}

At this point, we will structure the presentation in three parts regarding the Carleman estimate method and some trace operator inequalities.

\subsection{A first Carleman estimate}
	
	We start this discussion presenting a key Carleman inequality for solutions of 
\begin{equation}\label{adj-jlr}
\left\{\begin{array}{ll}
v_t+\left(x^\alpha v_x \right)_x=F, & (t,x)\in Q, \\
v(t,1)=0, & t\in (0,T),\\
\begin{cases*}
v(t,0)=0, & \text{ (Weak) } \\
\text{or}\\
({x^\alpha v_x})(t,0)=0,& \text{ (Strong) } 
\end{cases*},&  t\in (0,T), \\
v(T,x)=v_T(x), & x\in (0,1),
\end{array}\right.
\end{equation}
where $F\in L^2(Q)$ and $v_T\in L^2(0,1)$. We recall that \eqref{adj-jlr} is the adjoint system associated to \eqref{dcp}. Let us set $\psi :[0,1]\longrightarrow \R$, $\theta : (0,T) \longrightarrow \R$ and $\varphi :(0,T)\times [0,1] \longrightarrow \R$, given by 
\begin{equation}\label{functions}
\psi(x)=\frac{x^{2-\alpha}-3}{2-\alpha}, \ \theta(t):=\frac{1}{[t(T-t)]^4}, 
\varphi(t,x):=\theta(t)\psi(x),
\end{equation}
for any $(t,x)\in (0,T)\times [0,1]$. We observe that $\psi$, $\theta$ and $\varphi$ compose the weight functions which appear in \eqref{carl-frnt-inq}. Before, we state a Hardy-Poincaré type inequality, proved in \cite{alabau2006carleman} , which will be a very important ingredient from now on:

\begin{prop}
Assume $\alpha\in (0,2)$ and $\alpha\neq 1$. Let $w:[0,1]\longrightarrow \R$ be a locally absolutely continuous in $(0,1]$, with
\[\into x^\alpha |w_x|^2\,dx<+\infty.\]
Then, the following inequality holds
\begin{equation}\label{HP}
    \int_0^1x^{\alpha-2}|w|^2\,dx\leq \frac{4}{(1-\alpha)^2}\into x^\alpha |w_x|^2\,dx, 
\end{equation}
provided that either $\alpha\in (0,1)$ and $\dps\lim_{\ x\to 0^+}w(x)=0$ or $\alpha \in (1,2)$ and $\dps\lim_{\ x\to 1^-}w(x)=0$.
\end{prop}

\begin{prop}[Carleman Inequality]\label{car-frnt}
	There exist $C>0$ and $s_0>0$ such that, if $s\geq s_{0}$, then every solution  $v$ of (\ref{adj-jlr}) satisfies
	\begin{multline}
	\intq e^{2s\varphi}\left[(s\theta)^{-1}(|v_t|^2+|(x^\alpha v_x)_x|^2)+(s\theta) x^\alpha |v_x|^2+(s\theta)^{3}x^{2-\alpha}v^2 \right]\, dxdt \\
	\leq C\left(\intq e^{2s\varphi}|F|^2\,dxdt   +s\int_0^T e^{2s\varphi(t,1)}\theta(t) |v_x(t,1)|^{2}\,dt\right).\label{carl-frnt-inq}
	\end{multline}
\end{prop}

\begin{proof}
In \cite{alabau2006carleman}, the terms
\[
\intq e^{2s\varphi}(s\theta) x^\alpha |v_x|^2 dxdt
\emph{ and }
\intq e^{2s\varphi} (s\theta)^{3}x^{2-\alpha}v^2 dxdt
\]
have already been estimated as in \eqref{carl-frnt-inq}. 

In order to deal with the remaining terms, let us define 
$$w(t,x):=e^{s\varphi(t,x)}v(t,x),$$
where $v$ is the solution of (\ref{adj-jlr}). 
In Appendix  \ref{app}, we prove that there exists $C>0$ such that
\begin{equation}\label{uu}
\intq(s\theta)^{-1}(w_t^2+|(x^{\alpha} w_x)_x|^2)dxdt
\leq C\left(\intq e^{2s\varphi}|F|^2 dxdt+s \int_0^T \theta(t)w_x^2(t,1)dt \right)
\end{equation}
for any $s>0$ sufficiently large (see Lemma \ref{prop3.3}). For the next computations, observe that 
\[
\left\{ \begin{array}{l}
v_x= (-s\theta x^{1-\alpha} w+w_x)e^{-s\varphi}, \\
v_t=(-s\theta_t\psi w+w_t)e^{-s\varphi}.
\end{array}    \right.
\]

Firstly, since $\theta^{-1}\theta_t^2\psi^2\leq C\theta^{2}$, we apply \eqref{myeq12} to get
\begin{align}\label{myeq22}
&\intq e^{2s\varphi}(s\theta)^{-1}v_t^2 dxdt
=\intq (s\theta)^{-1}(-s\theta_t\psi w+w_t)^2 dxdt  \nonumber\\
& =\intq (s\theta)^{-1}s^2\theta_t^2\psi^2 w^2 dxdt  
+\intq (s\theta)^{-1}w_t^2 dxdt
-2\intq (s\theta)^{-1}s\theta_t\psi ww_t dxdt  \nonumber\\
& \leq C\intq (s\theta)^2 w^2 dxdt  
+\intq (s\theta)^{-1}w_t^2 dxdt
-2\intq \theta^{-1}\theta_t\psi ww_t dxdt \nonumber\\
&\leq \begin{multlined}[t]
C\left(\intq(s\theta)x^{\alpha} w_x^2 dxdt +\intq(s\theta)^{3} x^{2-\alpha} w^2 dxdt \right)\\ 
+\intq (s\theta)^{-1}w_t^2 dxdt 
-2\intq \theta^{-1}\theta_t\psi ww_t dxdt. 
\end{multlined}
\end{align}
Hence, using \eqref{uu}, \eqref{myeq12} and the estimate
\begin{align*}
& \left|\intq \theta^{-1}\theta_t\psi ww_t dxdt  \right|
\leq C\intq \theta^{-1/4}|\omega\omega_t| dxdt\\
& \leq C\intq |\theta^{-3/4}s^{-1/2}||\theta^{-1/2}s^{-1/2}w_t||s\theta w| dxdt \\
&\leq C\left( \intq (s\theta )^{-1}|w_t|^2 dxdt
+\intq (s\theta)^{2}|w|^2 dxdt \right) \\
& \leq  C\left( \intq (s\theta )^{-1}|w_t|^2 dxdt
+\intq(s\theta) x^{\alpha} |w_x|^2 dxdt +\intq (s\theta)^{3} x^{2-\alpha} |w|^2 dxdt \right),
\end{align*}
it is clear that 
\begin{align}\label{vt}
\intq e^{2s\varphi}(s\theta)^{-1}v_t^2 dxdt  & \leq C\left(\intq e^{2s\varphi}|F|^2 dxdt +s\int_0^T \theta(t)w_x^2(t,1)dt \right) \nonumber \\
&=C\left(\intq e^{2s\varphi}|F|^2 dxdt
+s\int_0^T \theta(t)v_x^2(t,1)e^{2s\varphi(t,1)} dt    \right),
\end{align}
recalling that 
\[
\displaystyle w_x(t,1)=(s\varphi_xe^{s\varphi}v+v_xe^{s\varphi})\vert_{x=1} =(s\varphi_xw+v_xe^{s\varphi})\vert_{x=1}=v_x(t,1)e^{s\varphi(t,1)}.
\]

Secondly, before dealing with 
$\intq e^{2s\varphi} (s\theta)^{-1} |(x^\alpha v_x)_x|^2dxdt,
$
we notice that
\[
 (x^{\alpha} v_x)_x = (-s\theta w-s\theta xw_x+(x^{\alpha} w_x)_x+(s\theta)^2 x^{2-\alpha} w-s\theta w_x x)e^{-s\varphi} .
\]
Thus, from \eqref{myeq12}, \eqref{uu} and Lemma \ref{prop3.3}, we readily get
\begin{align*}
& \intq (s\theta)^{-1}e^{2s\varphi}|(x^{\alpha} v_x)_x|^2 dxdt 
= \intq (s\theta)^{-1}|e^{s\varphi}(x^{\alpha} v_x)_x|^2  dxdt\nonumber\\
& \quad \leq \begin{multlined}[t]
    C\bigg(\intq (s\theta)^{-1}s^2\theta^2w^2 dxdt  
\quad +\intq (s\theta)^{-1}s^2\theta^2x^2w_x^2 dxdt  \\
 +\intq (s\theta)^{-1}|(x^{\alpha} w_x)_x|^2 dxdt  
+\intq (s\theta)^{-1}(s\theta)^4 x^{4-2\alpha} w^2 dxdt \\
 +\intq (s\theta)^{-1} (s\theta)^2 w_x^2 x^2 dxdt  \bigg)
\end{multlined}\\
&\quad \leq 
\begin{multlined}[t]
C\bigg(\intq (s\theta)^2 w^2 dxdt  
+\intq (s\theta) x^2 w_x^2 dxdt  \nonumber\\
+\intq (s\theta)^{-1}|(x^{\alpha} w_x)_x|^2 dxdt 
+\intq (s\theta)^{3} x^{2-\alpha} w^2 dxdt  
+\intq s\theta x^2 w_x^2 dxdt  \bigg)
\end{multlined}\\
&\quad\leq C\left(\intq e^{2s\varphi}|F|^2 dxdt +s\int_0^T \theta(t)w_x^2(t,1)dt \right) \nonumber \\
&\quad =C\left(\intq e^{2s\varphi}|F|^2 dxdt
+s\int_0^T \theta(t)v_x^2(t,1)e^{2s\varphi(t,1)} dt    \right).
\end{align*}
It means that \ref{carl-frnt-inq} holds.

\end{proof}

\subsection{Trace}\label{trace}

Given $a\in (0,1)$, we recall that $H^1(a,1)\hookrightarrow C^0([a,1])$ and 
\begin{equation}\label{TFC}
\displaystyle u(t)-u(s)=\int_s^tu'(\xi)\,d\xi,\     
\end{equation}
where $u\in H^{1} (a,1)$ and $s,t\in [a,1]$ (see \cite[Theorem 8.2]{Brezis}).
Let us set $\gamma_0:H^1(a,1)\longrightarrow \R$ by $\gamma_0 (u)=u(1)$, for each $u\in H^{1} (a,1)$. Clearly $\gamma_{0}$ is linear. Additionally, given $u,v\in H^1(a,1)$, we use \eqref{TFC} in order to obtain
\begin{align*}
|\gamma_{0} (u)-\gamma_{0} (v)|
&=|u(1)-v(1)| \\
&\leq |u(y)-v(y)|+\int_y^1|u' (\xi)-v'(\xi)|\,d\xi\\
& \leq |u(y)-v(y)|+(1-a)^{1/2} \n{u'-v'}{L^2(a,1)},
\end{align*}
for each $y\in [a,1)$. As a consequence,
\begin{align*}
    (1-a)|\gamma_{0} (u) - \gamma_{0} (v)|
    &\leq (1-a)^{1/2} \|u-v\|_{L^2 (a,1)}
    +(1-a)^{3/2} \|u^{\prime}-v^{\prime} \|_{L^2 (a,1)}
\end{align*}
implies
\begin{align*}
    |\gamma_{0} (u) - \gamma_{0} (v)|
    &\leq (1-a)^{-1/2} \|u-v\|_{L^2 (a,1)}
    +(1-a)^{1/2} \|u^{\prime}-v^{\prime} \|_{L^2 (a,1)}.
\end{align*}
In particular, 
\begin{align}\label{eqA1}
    |\gamma_{0} (u)|
    &\leq (1-a)^{-1/2} \|u \|_{L^2 (a,1)}
    +(1-a)^{1/2} \|u^{\prime} \|_{L^2 (a,1)} \nonumber \\
    &\leq [(1-a)^{-1/2} +(1-a)^{3/2}] \|u\|_{H^{1} (a,1)} 
\end{align}
for each $u\in H^{1} (a,1)$. It means that $\gamma_{0}$ is also continuous.

\begin{prop}\label{abalfa}
\begin{itemize}
    \item[(a)] There exists $A_{\alpha}>0$ such that 
    \[
    |u(1)|\leq A_{\alpha} \n{u}{H^1_\alpha}
    \]
for any $u\in H^1_\alpha$; 
\item[(b)] There exists $B_{\alpha}>0$ such that 
    \[
    |u_x (1)|\leq B_{\alpha} \n{u}{H^2_\alpha}
    \]
for any $u\in H^2_\alpha$.
\end{itemize}
\end{prop}

\begin{proof}
\begin{itemize}
\item[(a)]
Let us take $u\in H_{\alpha}^{2}$. For each $a\in (0,1)$, we notice that
\[u^{\prime} =\frac{1}{x^{\alpha/2}}x^{\alpha/2}u^{\prime}\leq \frac{1}{a^{\alpha/2}}x^{\alpha/2}u^{\prime} , \  \]
for any $x\in [a,1]$. In this case,
\begin{equation}\label{est.ux}
\n{u^{\prime}}{L^2(a,1)}^2\leq \frac{1}{a^{\alpha}}\n{x^{\alpha/2}u^{\prime}}{L^2(a,1)}^2\leq  \frac{1}{a^{\alpha}}\n{u}{H^1_\alpha}^2 
\end{equation}
and, consequently, we can use \eqref{eqA1} in order to get
\begin{align*}
    |u(1)|
    &\leq (1-a)^{-1/2} \|u \|_{L^2 (a,1)}
    +(1-a)^{1/2} \|u^{\prime} \|_{L^2 (a,1)} \\
    &\leq \left[ \frac{1}{(1-a)^{1/2}} +\frac{(1-a)^{1/2}}{a^{\alpha /2}} \right] \|u\|_{H_{\alpha}^{1}}.
\end{align*}
Since 
\[
\displaystyle f: \lambda \in (0,1) \longmapsto \left[ \frac{1}{(1-\lambda )^{1/2}} +\frac{(1-\lambda )^{1/2}}{\lambda ^{\alpha /2}} \right] \in (0,+\infty)
\]
is a continuous function satisfying $\displaystyle\lim_{\lambda \to 0^{+}} f(\lambda) = +\infty$ and $\displaystyle \lim_{\lambda \to 1^{-}} f(\lambda) = +\infty$,
there exists $\textbf{a} \in (0,1)$ such that $\displaystyle 0<f(\textbf{a} ) \leq \min_{\lambda (0,1)} f(\lambda)$. In particular, taking $A_{\alpha} = f(\textbf{a})$, the desired inequality follows.    

\vspace{0.5cm}

\item[(b)] For the second statement, we firstly fix $a\in (0,1]$ and observe that
\begin{equation*}
   u^{\prime \prime} =\frac{1}{x^\alpha}(x^\alpha u^{\prime})^{\prime}-\frac{\alpha}{x} u^{\prime} 
\end{equation*}
implies $|u^{\prime \prime}|^2\leq 4\left(\frac{1}{x^{2\alpha}}|(x^\alpha u^{\prime})^{\prime}|^2+\frac{\alpha^2}{x^2}|u^{\prime}|^2\right).$
Hence, using \eqref{est.ux}, we get
\begin{align}
\label{est.uxx}
\|u^{\prime \prime} \|_{L^2 (a,1)}^{2}
&\leq \frac{4}{a^{2\alpha}} \|u\|_{H_{\alpha}^{2}}^{2}
+\frac{4\alpha ^2}{a^2} \|u^{\prime} \|_{L^2 (a,1)}^{2} \nonumber \\
&\leq \left( \frac{4}{a^{2\alpha}} 
+\frac{4\alpha ^2}{a^{2+\alpha}} \right) \|u \|_{H_{\alpha}^{2}}^{2} 
\end{align}

As a consequence, \eqref{eqA1} and \eqref{est.ux} allow us to obtain
\begin{align*}
    \displaystyle |u^{\prime} (1)|
    &\leq (1-a)^{-1/2} \|u^{\prime} \|_{L^{2} (a,1)}
    +(1-a)^{1/2} \|u^{\prime \prime} \|_{L^{2} (a,1)} \\
    &\leq \left[ \frac{1}{a^{\alpha /2} (1-a)^{1/2}} + (1-a)^{1/2} \left( \frac{4}{a^{2\alpha}} + \frac{4\alpha ^2}{a^{2+\alpha}} \right)^{1/2}\right] \|u\|_{H_{\alpha}^{2}}.
\end{align*}
Arguing as before, we can take 
\[
\displaystyle B_{\alpha} = \min \left\{ \frac{1}{\lambda^{\alpha /2} (1-\lambda)^{1/2}} + (1-\lambda)^{1/2} \left( \frac{4}{\lambda^{2\alpha}} + \frac{4\alpha ^2}{\lambda^{2+\alpha}} \right)^{1/2}; \lambda \in (0,1) \right\} >0
\]
in order to complete the proof.

\end{itemize}
\end{proof}

Next, we deduce some consequences from the previous result:

\begin{cor}\label{abalfavec}
Let $A_{\alpha}$ and $B_{\alpha}$ be those two positive constants obtained in Proposition \ref{abalfa}. Then:
\begin{itemize}
    \item[(a)]  
    $\n{u(t,1)}{L^2(0,T)}\leq A_{\alpha} \n{u}{L^2(0,T;H^1_\alpha)}$
for any $u\in L^2(0,T; H^1_\alpha)$; 
\item[(b)]
$
\n{u_x(t,1)}{L^2(0,T)}\leq B_{\alpha} \n{u}{L^2(0,T;H^2_\alpha)}
$
for any $u\in L^2(0,T;H^2_\alpha)$.
\end{itemize}
\end{cor}

\begin{prop}\label{weak.conv}
Let $(u_n)_{n=1}^{\infty}$ be a sequence in $L^2(0,T;H^2_\alpha)$ which weakly converges to $u\in L^2(0,T;H^2_\alpha)$.
Then
\begin{equation*}
         u_n(t,1)\rightharpoonup u(t,1) \text{ and } u_{nx}(t,1)\rightharpoonup u_x(t,1)
\end{equation*}
weakly in $L^2 (0,T)$.
\end{prop}

\begin{proof}
Let us set $\Gamma_{i} : L^2 (0,T;H_{\alpha}^{2}) \longrightarrow L^2 (0,T)$, with $i\in \{1,2\}$, given by 
\[
\Gamma_1 (u)=u(\cdot,1) \emph{  and  } \Gamma_2 (u)= u_x (\cdot , 1), 
\]
for each $u\in  L^2 (0,T;H_{\alpha}^{2})$. Notice that Proposition \ref{abalfavec} means that $\Gamma_1$ and $\Gamma_2$ are two continuous linear mappings. In this case, for each $S\in (L^2 (0,T))^{\prime}$, $S\circ \Gamma_1$ and $S\circ \Gamma_2$ belong to $(L^2 (0,T;H_{\alpha}^{2}))^{\prime}$. Hence,
\[
S(u_n (\cdot ,1)) = (S\circ \Gamma_1)(u_n) \to (S\circ \Gamma_1)(u) = S(u (\cdot ,1))
\]
and
\[
S(u_{nx} (\cdot ,1)) = (S\circ \Gamma_2)(u_n) \to (S\circ \Gamma_2)(u) = S(u_x (\cdot ,1))
\]
as $n\to \infty$. It ends the proof.
\end{proof}

\subsection{Some useful technical results}
Let us obtain two identities and a convergence result that will be crucial in Section \ref{principal}.  Given $a\in (0,1)$, we recall the continuous embedding $H^1(a,1)\hookrightarrow C^0([a,1])$. Thus, for each $v\in L^2 (0,T;H_{\alpha}^{2})$, we can apply \eqref{TFC} to obtain
\begin{equation}\label{vx}
v_x(t,s)=v_x(t,1)-\int_s^1{ v_{xx}(t,r)} dr,    
\end{equation}
where $s\in [a,1]$ and $t\in (0,T)$. In particular, note that $v_x\in L^2(0,T;C^0([a,1])$. Integrating from $x\in [a,s]$ to 1, with respect to $s$, we readily get
\begin{equation}\label{eq.v}
    v(t,x)=-(1-x)v_x(t,1)+\int_x^1\int_s^1 {v_{xx}(t,r)}\,drds.
\end{equation}
for any $(t,x)\in (0,T)\times [a,1)$. 

Now, define
\[
V(t,x)=\frac{1}{1-x}\int_x^1\int_s^1v_{xx}(t,r)\,drds.
\]
Let us prove that 
\begin{equation}\label{lim-sup}
\sup_{x\in[a,1]}\n{V(\cdot,x)}{L^2(0,T)}^2\to 0, \text{ as } a\to 1^+.
\end{equation}
Indeed, 
from \eqref{vx}, define 
\[
f(t,s):=\int_s^1v_{xx}(t,r)dr=v_x(t,1)-v_x(t,s)\in L^2(0,T;C^0([a,1])).
\]

From Jensen's inequality, we take
\begin{equation*}
    |V(t,x)|^2 = \frac{1}{1-x}\int_x^1|f(t,s)|^2\, ds\leq \n{f(t,\cdot)}{C^0([a,1])}.
\end{equation*}
Integrating over $[0,T]$, we have that $\n{V(\cdot,x)}{L^2(0,T)}^2\leq \n{f}{L^2(0,T;C^0([a,1]))}^2,\ \forall x\in [a,1]$, whence 
\[
\dps\sup_{x\in[a,1]}\n{V(\cdot,x)}{L^2(0,T)}^2\leq  \n{f}{L^2(0,T;C^0([a,1]))}^2.
\]
Therefore, we just need to prove that $\n{f}{L^2(0,T;C^0([a,1]))}^2\to 0$, as $a\to 1^-.$ To do this, we need the following lemma.

\begin{lem}\label{lema-sup}
Given $w \in C^0([a,1])$, define $g(\xi)=\dps\sup_{s\in[\xi,1]}|w(1)-w(s)|$ for each $\xi \in [a,1]$. Then, $\dps\lim_{\xi\to 1^-}g(\xi)=0$.
\end{lem}
\begin{proof}
Since $w$ is continuous in $s=1$, given $\ep>0$, there exists $\delta>0$ such that $|w(1)-w(s)|<\ep/2$, $\forall s\in (1-\delta,1)$. In particular, if $\xi\in (1-\delta,1)$, we have that
\[|w(1)-w(s)|<\ep/2,\ \forall s\in (\xi,1),\]
whence $g(\xi)=\dps\sup_{s\in[\xi,1]}|w(1)-w(s)|\leq \ep/2<\ep$.
\end{proof}
Applying Lemma \ref{lema-sup} for $w=v_x(t,\cdot)$, we have
\[
\n{f(t,\cdot)}{C^0([a,1])}=\sup_{s\in [a,1]}|v_x(t,1)-v_x(t,s)|\to 0, \text{ as } a\to 1^-,\ \text{ a.e in }  [0,T].
\]
Since $\n{f(t,\cdot)}{C^0([a,1])}\leq 2\n{v_x(t,\cdot)}{C^0([a,1)} \leq C\n{v_x(t,\cdot)}{H^{1} (a,1)}\in L^2(0,T)$, Lebesgue's dominated convergence theorem gives us \eqref{lim-sup}.

\subsection{A refined Carleman estimate}

At this place, following the ideas presented in \cite{chaves2020boundary}, we will improve that Carleman estimate proved in Proposition \ref{car-frnt} by emphasizing the influence of the control domain in \eqref{dcp}. Precisely, we aim the following result:

 \begin{thm}
     There exist positive constants $C>0$ and $s_0>0$, only depending on $T$, such that, if $s\geq s_0$, then every solution $v$ of \eqref{adj-jlr} verifies 
     \begin{multline}
 	\intq e^{2s\varphi}\left[(s\theta)^{-1}(|v_t|^2+|(x^\alpha v_x)_x|^2)+(s\theta) x^\alpha |v_x|^2+(s\theta)^{3}x^{2-\alpha}v^2 \right]\, dxdt \\
 	\leq C\left(\intq e^{2s\varphi}|F|^2\,dxdt   +\frac{1}{\ep^3}s^7\intw e^{2s\varphi}\theta^7 x^{2-\alpha}|v|^{2}\,dxdt\right).\label{carl-dist}
 	\end{multline}
 \end{thm}

 	\begin{proof}
 	Comparing \eqref{carl-dist} to \eqref{carl-frnt-inq}, the only difference appears in the last term of the right side. So that, this current proof is concerned the obtainment of
\[
 \frac{1}{\ep^3}s^7\intw e^{2s\varphi}\theta^7 x^{2-\alpha}|v|^{2}\,dxdt
\]
from
\[
\displaystyle s\int_0^T e^{2s\varphi(t,1)}\theta(t) |v_x(t,1)|^{2}\,dt,
\]
being aware that each extra term arised from the estimations must be incorporated to that ones already given in \eqref{carl-frnt-inq}. 
In order to simplify the notation, let us define the right hand side of \eqref{carl-frnt-inq} as 
\[I(s,v)=\intq e^{2s\varphi}\left[(s\theta)^{-1}(|v_t|^2+|(x^\alpha \xi_x)_x|^2)+(s\theta) x^\alpha |v_x|^2+(s\theta)^{3}x^{2-\alpha}v^2 \right]\, dxdt,\]
for any $v$  solution \eqref{adj-jlr}.

Firstly, as in the proof of Proposition \ref{car-frnt}, by a standard density argument, we can consider $v$ as a solution of \eqref{adj-jlr} which is sufficiently regular. Now, let us take a cut-off function $\zeta\in C^3([0,1])$ such that $\zeta=1$ in $(1-\frac{\varepsilon}{2},1)$, $\zeta=0$ in $[0,1-\varepsilon]$, with $$\zeta'=O(\varepsilon^{-1}), \ \ \ \zeta''=O(\varepsilon^{-2}) \ \ \ \mbox{and} \ \ \ \zeta'''=O(\varepsilon^{-3}).$$ Thus, a simple computation gives
 \begin{align*}
s\theta(t)e^{2s\varphi(t,1)}|v_x(t,1)|^2
&=s\int_{1-\varepsilon}^1\displaystyle\frac{d}{dx}[\theta(t)e^{2s\varphi(t,x)}x^{2\alpha}\zeta(x)|v_x(t,x)|^2]\,dx \\
& =\begin{multlined}[t]
2s^2\int_{1-\varepsilon}^1\theta(t)^2e^{2s\varphi(t,x)} x^{{\alpha+1}}\zeta(x)|v_x(t,x)|^2\,dx \\  
+s\int_{1-\varepsilon}^1\theta(t)e^{2s\varphi(t,x)}x^{2\alpha}\zeta'(x)|v_x(t,x)|^2\,dx\\
 +2s\int_{1-\varepsilon}^1\theta(t)e^{2s\varphi(t,x)}\zeta(x)x^\alpha v_x(t,x)(x^\alpha v_x(t,x))_x\,dx.
\end{multlined}
\end{align*}
Integrating over $[0,T]$, we obtain
\begin{align*}
     s\int_0^Te^{2s\varphi(t,1)}\theta(t)|v_x(t,1)|^2\, dt 
=&\begin{multlined}[t]
2s^2\int_0^T\int_{1-\varepsilon}^1\theta^2e^{2s\varphi} x^{\alpha+1}\zeta|v_x|^2\,dx dt 
  +s\int_0^T\int_{1-\varepsilon}^1\theta e^{2s\varphi}x^{2\alpha}\zeta'|v_x|^2\,dxdt\\
+2s\int_0^T\int_{1-\varepsilon}^1\theta e^{2s\varphi}\zeta x^\alpha v_x(x^\alpha v_x)_x\,dxdt
\end{multlined}\\
& =:I_1+I_2+I_3.
 \end{align*}
 
As aforementioned, it suffices to bound each integral by terms of the form
\[\delta I(s,\varphi) + \frac{C}{\delta}\frac{s^7}{\ep^3}\intw e^{2s\varphi}\theta^7 x^{2-\alpha}|v|^{2}\,dxdt,\]
 where $\delta>0$ must be chosen sufficiently small.
 
 Indeed, since $\theta^2\leq C\theta^3$ and $s^2\leq Cs^3$, then we have
 \begin{equation*}
    I_1  \leq Cs^3\intw e^{2s\varphi}\theta^3x^\alpha |v_x|^2\zeta\,dxdt.\\
     \end{equation*} 

From Young's inequality with $\delta$, we also have
 \begin{align*}
    I_3 & =2\int_0^T\int_{1-\varepsilon}^1\left( e^{s\varphi}(s\theta)^{-1/2}  (x^\alpha v_x)_x\zeta^{1/2}\right)\left(e^{s\varphi}(s\theta)^{3/2} x^{\alpha} v_x\zeta^{1/2} \right)\,dxdt\\
    &\leq \delta \intw e^{2s\varphi} (s\theta)^{-1}  |(x^\alpha v_x)_x|^2\zeta \,dxdt+ \frac{C}{\delta}s^3\int_0^T\int_{1-\varepsilon}^1 e^{2s\varphi}\theta^3 x^{2\alpha} |v_x|^2\zeta \\
    &\leq \delta I(s,\varphi)+\frac{C}{\delta}s^3\intw e^{2s\varphi}\theta^3x^\alpha|v_x|^2\zeta\,dxdt.
\end{align*} 

Hence,
\begin{equation*}
    I_1+I_3\leq \delta I(s,\varphi)+Cs^3\intw e^{2s\varphi}\theta^3x^\alpha|v_x|^2\zeta\,dxdt.
\end{equation*}
Now, integrating by part the last integral, we have
 \begin{align*}
 s^3\int_0^T\int_{1-\varepsilon}^1
 e^{2s\varphi}\theta^3\zeta x^\alpha|v_x|^2\,dx\,dt
&  =s^3\int_0^T\int_{1-\varepsilon}^1
 (e^{2s\varphi}\theta^3\zeta x^\alpha v_x)v_x\,dx\,dt\nonumber\\
 & =\begin{multlined}[t]
  -2s^4\int_0^T\int_{1-\varepsilon}^1e^{2s\varphi}\theta^4\zeta xvv_x\,dxdt\\  -s^3\int_0^T\int_{1-\varepsilon}^1e^{2s\varphi}\theta^3\zeta'x^\alpha vv_x\,dxdt\\
  -s^3\int_0^T\int_{1-\varepsilon}^1e^{2s\varphi}\theta^3\zeta v(x^\alpha v_x)_x\,dx\,dt
  \end{multlined}\nonumber\\ 
 &=: J_1+J_2+J_3.
 \end{align*}
Using Young's inequality with $\delta$ again, 
\begin{align*}
J_1&=\intw\left(e^{s\varphi}(s\theta)^{1/2}x^{\alpha/2}v_x\right) \left(2e^{s\varphi}(s\theta)^{7/2}x^{1-\alpha/2}v \zeta\right)\,dxdt\\
&\leq \delta\int_0^T\int_{1-\varepsilon}^1 e^{2s\varphi}(s\theta) x^\alpha|v_x|^2\,dx\,dt+\frac{C}{\delta}s^7\int_0^T\int_{1-\varepsilon}^1e^{2s\varphi}\theta^7x^{2-\alpha}|v|^2\zeta^2\,dx\,dt\\
&\leq \delta I(s,\varphi)+\frac{C}{\delta}s^7\int_0^T\int_{1-\varepsilon}^1e^{2s\varphi}\theta^7x^{2-\alpha}|v|^2\,dx\,dt,
\end{align*}
Similarly,
\begin{align*}
 J_2&\leq \delta\int_0^T\int_{1-\varepsilon}^1e^{2s\varphi}(s\theta) x^\alpha|v_x|^2\,dx\,dt+\frac{C}{\delta}s^5\varepsilon^{-2}\int_0^T\int_{1-\varepsilon}^1e^{2s\varphi}\theta^5|v|^2\,dx\,dt,\nonumber\\
 J_3&\leq \delta\int_0^T\int_{1-\varepsilon}^1e^{2s\varphi}(s\theta)^{-1}|(x^\alpha v_x)_x|^2\,dx\,dt+\frac{C}{\delta}s^7\int_0^T\int_{1-\varepsilon}^1e^{2s\varphi}\theta^7|v|^2\,dx\,dt.
\end{align*}
Since $1-\ep\leq x\leq 1$, we have $x^{2-\alpha}\geq C>0$, then
\begin{equation*}
    I_1+I_3\leq \delta I(s,\varphi)+\frac{C}{\delta}s^7\varepsilon^{-3}\int_0^T\int_{1-\varepsilon}^1e^{2s\varphi}\theta^7x^{2-\alpha}|v|^2\,dxdt.
\end{equation*}

It just remains to bound $I_2$. Integrating by parts, we see that
\begin{align*}
 I_2=& \begin{multlined}[t]
 -2s^2\int_0^T\int_{1-\varepsilon}^1e^{2s\varphi}\theta^2\zeta'x^{\alpha+1}vv_x\,dxdt
 -s\int_0^T\int_{1-\varepsilon}^1e^{2s\varphi}\theta\zeta''x^{2\alpha}vv_x\,dxdt\\
 -s\alpha\int_0^T\int_{1-\varepsilon}^1e^{2s\varphi}\theta\zeta'x^{2\alpha-1} vv_x\,dx\,dt
 -s\int_0^T\int_{1-\varepsilon}^1e^{2s\varphi}\theta\zeta'x^\alpha v(x^\alpha v_x)_x\,dxdt
 \end{multlined} \nonumber \\ 
 =:&  K_1+K_2+K_3+K_4
\end{align*}
As before,
\begin{equation*}
    K_1\leq \delta \int_0^T\int_{1-\varepsilon}^1e^{2s\varphi}(s\theta) x^\alpha|v_x|^2\,dx\,dt +\frac{C}{\delta}\varepsilon^{-2}s^3\int_0^T
 e^{2s\varphi}\theta^3x^{2-\alpha}|v|^2\,dx\,dt,
\end{equation*}
\begin{align*}
 K_2&= \frac{s}{2}\int_0^T\int_{1-\varepsilon}^1e^{2s\varphi}\theta\zeta''x^{2\alpha}\frac{d}{dx}(|v|^2)\,dx\,dt\\ 
  & =\begin{multlined}[t]
  -s^2\int_0^T\int_{1-\varepsilon}^1e^{2s\varphi}\theta^2x^{\alpha+1}\zeta''|v|^2\,dx\,dt-\frac{s}{2}\int_0^T\int_{1-\varepsilon}^1e^{2s\varphi}\theta\zeta'''x^{2\alpha}|v|^2\,dx\,dt\\
  -\alpha s\int_0^T\int_{1-\varepsilon}^1
  e^{2s\varphi}\theta\zeta''x^{2\alpha-1}|v|^2\,dx\,dt
  \end{multlined}\\
 &\leq Cs^2\varepsilon^{-3}\int_0^T\int_{1-\varepsilon}^1
 e^{2s\varphi}\theta^2|v|^2\,dx\,dt,
  \end{align*}
 \begin{equation*}
 K_3\leq \delta\int_0^T\int_{1-\varepsilon}^1
 e^{2s\varphi}(s\theta) x^\alpha|v_x|^2\,dx\,dt+\frac{C}{\delta}s\varepsilon^{-2}\int_0^T\int_{1-\varepsilon}^1 e^{2s\varphi}\theta|v|^2\,dx\,dt,
 \end{equation*}
 \begin{equation*}
 K_4\leq \delta\int_0^T\int_{1-\varepsilon}^1e^{2s\varphi}(s\theta)^{-1}|(x^\alpha v_x)_x|^2\,dx\,dt+\frac{C}{\delta}s^3\varepsilon^{-2}\int_0^T\int_{1-\varepsilon}^1e^{-2s\varphi}\theta^3|v|^2\,dx\,dt.
  \end{equation*}
 
 Therefore, we also have 
 \[I_2\leq \delta I(s,\varphi)+\frac{C}{\delta}s^7\varepsilon^{-3}\int_0^T\int_{1-\varepsilon}^1e^{2s\varphi}\theta^7x^{2-\alpha}|v|^2\,dxdt,\]
 as required.
 \end{proof}

\section{Internal null Controllability for the degenerate problem}\label{familia}

The aim of this section is to obtain a family $(\yep,\fep)_{\ep >0}$ of solutions to the internal null  problem \eqref{dcp}. In order to do that, we will first define some weights that does not vanish at $t=0$, using them to prove a certain observability inequality. It will allow us to solve the control problem \eqref{dcp}.

Let us consider a function $m\in C^\infty([0,T])$ satisfying
\[\left\{\begin{array}{ll}
m(t)\geq t^4(T-t)^4, & t\in (0, T/2];\\
m(t)= t^4(T-t)^4, & t\in \left[T/2,T\right];\\
m(0)>  0, &
\end{array}
\right.\]
and define	
\begin{equation}\label{functions2} 
\begin{gathered}
        \tau(t):=\frac{1}{m(t)},\ 
A(t,x):=\tau (t)\psi(x),\\
A^\ast(t):=\max_{0\leq x\leq 1}A(t,x), \
\hat{A}(t,x) :=\min_{0\leq x\leq 1}A(t,x) \mbox{ and } \bar{A}(t):=2A^\ast-\hat{A},
\end{gathered}
\end{equation}
where $(t,x)\in [0,T)\times [0,1]$. Observe that $\bar{A}<0$. Under these definitions and notations, we have the next result: 

\begin{lem}\label{lem3.1}
There exists $C>0$ such that
\begin{multline}\label{car-3.1}
    \n{v(0)}{L^2(0,1)}^2+\intq e^{2sA}\bigg[\tau^{-1} \left(|v_t|^2+|(x^\alpha v_x)_x|^2\right)+\tau x^\alpha |v_x|^2+\tau^3x^{2-\alpha}|v|^2 \bigg]\,dxdt\\
    \leq C\left(\ep^{-3}\intw e^{2sA}\tau^7x^{2-\alpha}|v|^2\,dxdt+\intq e^{2sA}|F|^2\,dxdt\right),
\end{multline}
for every $v$ solution of \eqref{adj-jlr}.
\end{lem}

\begin{proof} Firstly, since $\varphi\leq A$, $\varphi=A$ in $[T/2,T]$  and $e^{2sA}\tau^7\geq C>0$ in $[0,T/2]$, we notice that
\begin{multline}\label{eq.4.3}
 \intq e^{2s\varphi}|F|^2\,dxdt   +\frac{1}{\ep^3}s^7\intw e^{2s\varphi}\theta^7 x^{2-\alpha}|v|^{2}\,dxdt\\
      \leq \intq e^{2sA}|F|^2\,dxdt   +\frac{1}{\ep^3}s^7\intw e^{2sA}\tau^7 x^{2-\alpha}|v|^{2}\,dxdt.
\end{multline}

 As usual, we will divide the proof in two  different situations: when $t\in [0,T/2]$ and when  $t\in [T/2,T]$. The second case follows immediately from Carleman inequality \eqref{carl-dist} and the last inequality \eqref{eq.4.3},  since $A=\varphi$ and $\theta=\tau$ in $[T/2,T]$.
 

In order to obtain the inequality in $[0,T/2]$, let us consider a cut-off function $\xi\in C^\infty([0,T]$ such that
\[0\leq \xi\leq 1,\ \xi=1 \in [0,T/2] \text{ and } \xi=0 \in [3T/4,T].\]
Define $w=\xi v$ and note that $w$ is a solution to the problem
\begin{equation*}
    \begin{cases}
    -w_t-(x^\alpha w_x)_x=\xi F-\xi'v,  & (t,x) \in Q,\\
    w(t,1)=0, & \in t\in  (0,T),\\
    \begin{cases}
    w(t,0)=0\\
    \text{or}\\
    (x^\alpha w_x)(t,0)
    \end{cases} ,& t \in (0,T),\\
    w(T,x)=0 ,& x\in (0,1).
    \end{cases}
\end{equation*}
From Proposition \ref{well-pos}, 

\[w\in H^1(0,T;L^2(0,1))\cap L^2(0,T;H^2_\alpha)\cap C^0([0,T];H^1_\alpha),\]
and
\[\sup_{0\leq t\leq T} \n{w(t)}{H^1_\alpha}^2+\n{w_t}{L^2(Q)}^2+\n{(x^\alpha w_x)_x}{L^2(Q)}^2\leq C\left(\n{\xi F-\xi'v}{L^2(Q)}^2\right).\]
As a consequence, 
\[
\n{w(0)}{L^2(0,1)}^2 \leq C\left(\n{\xi F-\xi'v}{L^2(Q)}^2\right)
\]
and
\[
 \n{w}{L^2(0,T;H^2_\alpha)}^2\leq C\left(\n{\xi F-\xi'v}{L^2(Q)}^2\right).
\]
Hence, since $v=w$ in [0,T/2], we have
\begin{equation*}
    \n{v(0)}{L^2(0,1)}^2+\n{v}{L^2(0,T/2;H^2_\alpha)}^2+\n{v_t}{L^2(0,T/2;L^2(0,1))}^2\leq C\left(\n{\xi F-\xi'v}{L^2(Q)}^2\right).
\end{equation*}
Now, we will estimate the right side of this last inequality. Recalling that $e^{2sA}$ is bounded from below in $[0,T/2]$, we obtain 
\[\n{\xi F}{L^2(Q)}^2 \leq \int_0^{T/2}\int_0^1 e^{2sA}|F|^2\,dxdt\leq C\left(\ep^{-3}\intw e^{2sA}\tau^7x^{2-\alpha}|v|^2\,dxdt+\intq e^{2sA}|F|^2\,dxdt\right).\]
Besides, since $e^{2s\varphi}(s\theta)^k$ is bounded from below in $[T/2, 3T/4]$, for $k\in \mathbb{Z}$, we can use \eqref{HP}, \eqref{carl-dist} and \eqref{eq.4.3}, achieving

\begin{align*}
    \n{\xi'v}{L^2(Q)}^2 & \leq \int_{T/2}^{3T/4}\int_0^1 |v|^2\, dxdt\\
                        & \leq C\left(\int_{T/2}^{3T/4}\int_0^1e^{2s\varphi} (s\theta)^3x^{2-\alpha}|v|^2\,dxdt+\int_{T/2}^{3T/4}\int_0^1e^{2s\varphi} (s\theta)x^{\alpha}|v_x|^2\,dxdt\right)\\
                        & \leq C\left(\ep^{-3}\intw e^{2sA}\tau^7x^{2-\alpha}|v|^2\,dxdt+\intq e^{2sA}|F|^2\,dxdt\right).
\end{align*}
Finally, the boundedness of $e^{2sA}\tau^k$, for $k\in \mathbb{Z}$, in $[0,T/2]$, implies
\begin{align*}
      & \n{v(0)}{L^2(0,1)}^2+\int_0^{T/2}\int_0^1 e^{2sA}\bigg[\tau^{-1} \left(|v_t|^2+|(x^\alpha v_x)_x|^2\right)+\tau x^\alpha |v_x|^2+\tau^3x^{2-\alpha}|v|^2 \bigg]\,dxdt\\
    & \hspace{1cm}\leq \n{v(0)}{L^2(0,1)}^2+\n{v}{L^2(0,T/2;H^2_\alpha)}^2+\n{v_t}{L^2(0,T/2;L^2(0,1))}^2\\ 
    &\hspace{1cm} \leq C\left(\ep^{-3}\intw e^{2sA}\tau^7x^{2-\alpha}|v|^2\,dxdt+\intq e^{2sA}|F|^2\,dxdt\right),
\end{align*}
which ends the proof.
\end{proof}
We finish this section constructing a family of solution to the distributed null controllability problem \eqref{dcp}. Initially, let us define  the operator $\mathcal{L}^\ast w=-w_t-(x^\alpha w_x)_x$,  the linear spaces 
\begin{equation*}
    \begin{cases}
    P_\alpha=\{w\in C^2(\bar{Q}); \ w(t,1)=0, \text{ and } w(t,0)=0,\ \forall t\in (0,T)\}, &  \text{ if } \alpha \in (0,1);\\[0.5em]
 P_\alpha= \{w\in C^2(\bar{Q}); \ w(t,1)=0, \text{ and } x^\alpha w(t,0)=0,\ \forall t\in (0,T)\}, &    \text{ if } \alpha \in [1,2),
    \end{cases}
\end{equation*}
and the bilinear form
    \begin{equation*}
    a_\ep(w_1,w_2)=\intq e^{2sA}\mathcal{L}^\ast w_1\mathcal{L}^\ast w_2\, dxdt+\frac{1}{\ep^3}\intw e^{2sA^\ast}\tau^7w_1w_2\, dxdt,
\end{equation*}
for each $\varepsilon \in (0,1)$. It is a consequence of Carleman inequality \eqref{car-3.1} that each $a_\ep$ is symmetric and positive definite. So that, let $P_{\varepsilon}$ be the completion of $P_{\alpha}$ with respect to the inner product $a_{\varepsilon}$, and denote by $\n{\cdot}{P_\ep}$ the corresponding norm. 
Taking the linear form
\[
\displaystyle \ell : \varphi \in P_{\ep} \longmapsto \int_{0}^{1} u_0 (x) \varphi (0,x)dx \in \R,
\]
we can apply Carleman inequality \eqref{car-3.1} to see that
\[
|\langle \ell, \varphi\rangle|\leq \n{u_0}{L^2(0,1)}\n{\varphi(0)}{L^2(0,1)}\leq  C\n{u_0}{\sqrt{a_\ep(\varphi,\varphi)}}=C\n{u_0}{L^2(0,1)}\n{\varphi}{P_\ep},
\]
where $C$ does not depend on $\ep$.
Therefore, by Lax-Milgram Theorem, there exists a unique $\hat{\varphi}_\ep\in P_\ep$ such that 
\[a_\ep(\hat{\varphi}_\ep,\varphi)=\langle \ell,\varphi\rangle,\ \forall \varphi\in P_\ep,\]
that is,
\[\intq e^{2sA}\mathcal{L}^\ast \vep\mathcal{L}^\ast \varphi\, dxdt +\intw\frac{1}{\ep^3}e^{2sA^\ast}\tau^7\vep \varphi\,dxdt=\into u_0(x)\varphi(0,x)\, dx,\ \forall \varphi \in P_\ep.\]
It means that, for each $\varepsilon \in (0,1)$, 
\begin{equation*}
\yep:=e^{2sA}\mathcal{L}^\ast \vep  \text{ and } \fep:=-\frac{1}{\ep^3}e^{2sA^\ast}\tau^7\vep
\end{equation*}
are a state and a control, respectively, which solve the distributed null controllability problem \eqref{dcp}.
Indeed, for any $(F,z^T)\in L^2(Q)\times L^2(0,1)$, if $z$ is a solution of \eqref{adj-transp}, we have $z\in P_\ep$ and then
\begin{equation}\label{sol.dcp}
\intq \yep F\, dxdt =\intw \fep z\,dxdt+\into u_0(x)z(0,x)\, dx.
\end{equation}
Furthermore, $\n{\vep}{P_\ep}\leq C\n{u_0}{L^2(0,1)}$, since
\begin{equation*}
    \n{\vep}{P_\ep}^2=a_\ep(\vep,\vep)=\langle \ell,\vep \rangle \leq C\n{u_0}{L^2(0,1)}\n{\vep}{P_\ep}.
\end{equation*}
Thus, we can use $\mathcal{L}^\ast\vep=e^{-2sA}\yep$ and $\vep=-\ep^3e^{-2sA^\ast}\tau^{-7}\fep$ to conclude that 
\begin{align*}
    \displaystyle 
    &\intq e^{-2sA}|\yep|^2\,dxdt + {\ep^3}\intw e^{-2sA^\ast}\tau^{-7}|\fep|^2\,dxdt
    =\n{\vep}{P_\ep}^2
    \leq C^2\n{u_0}{L^2(0,1)}^2.
\end{align*}
Therefore,
\begin{equation}\label{3.3*}
\begin{split}
        & \|e^{-sA^\ast}\tau^{-7/2}\fep\|_{L^2((0,T)\times \omega_\ep)}\leq {\frac{C}{\ep^3}}\n{u_0}{L^2(0,1)}.\\
        & \|e^{-sA}\yep\|_{L^2(Q)}\leq C\n{u_0}{L^2(0,1)},
\end{split}
\end{equation}
where $C$ does not depend on $\ep$. As a result, each $(\yep,\fep)$ is really a state-control solution to \eqref{dcp}.


\section{Proof of the main result}\label{principal}

In this section, in addition to prove the theorem, we make it clear how the convergence of the distributed control problem \eqref{dcp} to the boundary control problem \eqref{bcp} works. 

In the previous section we have obtained a family $(\yep,\fep)_{\ep>0}$ of solutions to problem \eqref{dcp} in the sense of Definition \ref{transp}, that is to say that $(\yep,\fep)$ satisfies \eqref{sol.dcp}. In order to obtain a solution to the boundary null control problem \eqref{bcp}, we need to pass limits in that equation and obtain $(u,h)$ that satisfies equation \eqref{sol.bound}. 

The convergence of the left hand side of \eqref{sol.dcp} is a immediate consequence of estimate \eqref{3.3*}. Next, we will prove two Lemmas that will guarantee the convergence of the right hand side. 

\begin{lem}\label{prop4.1}
For each $\varepsilon \in (0,1)$, $e^{s\hat{A}}\tau^{-1/2}\vep \in L^2(0,T;H^2_\alpha)$, with 
\begin{equation*}
    \n{e^{s\hat{A}}\tau^{-1/2}\vep}{L^2(0,T;H^2_\alpha)}\leq C\n{u_0}{L^2(0,1)},
\end{equation*}
where $C$ does not depend on $\ep$. Moreover, there exists a function $\varphi$ such that, up to a subsequence, 
\begin{equation*}
         e^{s\hat{A}}\tau^{-1/2}\vep \rightharpoonup e^{s\hat{A}}\tau^{-1/2}\varphi \text{ weakly in } L^2(0,T;H^2_\alpha).  \\
\end{equation*}
As a consequence,
\begin{equation*}
     e^{2sA^\ast}\tau^{7}\vepx(t,1) \rightharpoonup e^{2sA^\ast}\tau^{7}\varphi_x(t,1) \text{ weakly in } L^2(0,T).
\end{equation*}
\end{lem}

\begin{proof}
Firstly, since $\tau\geq C>0$, we obtain $e^{2s\hat{A}}\tau^{-1}\leq C e^{2sA}\tau^k$, $\forall k\geq -1$. From Hardy-Poincar\'e \eqref{HP} and  Lemma  \ref{lem3.1}, we have
\begin{align*}
    \n{e^{s\hat{A}}\tau^{-1/2}\vep}{L^2(0,T;H^2_\alpha)}^2& =\int_0^Te^{2s\hat{A}}\tau^{-1} \left(\n{\vep(t)}{L^2(0,1)}^2+\n{x^{\alpha/2}\vepx(t)}{L^2(0,1)}^2+\n{(x^\alpha \vepx)_x(t)}{L^2(0,1)}^2\right)dt\\
    & \leq C\intq e^{2sA}\left(\tau^{-1}|(x^\alpha \vepx)_x|^2+\tau x^\alpha |\vepx|^2+\tau^3x^{2-\alpha}|\vep|^2\right)\,dxdt\\
    & \leq C\left(\ep^{-3}\intw e^{2sA}\tau^7x^{2-\alpha}|\vep|^2\,dxdt+\intq e^{2sA}|\mathcal{L}\vep|^2\,dxdt\right)\\
    & \leq a(\vep,\vep)\leq C\n{u_0}{L^2(0,1)}^2.
\end{align*}
Therefore, up to a subsequence, we have
\begin{equation*}
    e^{s\hat{A}}\tau^{-1/2}\vep \rightharpoonup e^{s\hat{A}}\tau^{-1/2}\varphi \text{ weakly in } L^2(0,T;H^2_\alpha).
\end{equation*}
From Proposition \ref{weak.conv}, we get
\begin{equation}\label{wcx}
    e^{s\hat{A}}\tau^{-1/2}\vepx(t,1)\rightharpoonup e^{s\hat{A}}\tau^{-1/2}\varphi_x(t,1) \text{ in } L^2(0,T).
\end{equation}
In order to  prove the last convergence, we notice that $e^{2sA^\ast}\tau^7=e^{s\hat{A}}\tau^{-1/2}(e^{s\bar{A}}\tau^{15/2})=e^{s\hat{A}}\tau^{-1/2}w(t)$, where $w(t)$ is bounded. Thus, given $v\in L^2(0,T)$, we have $wv\in L^2(0,T)$. Therefore, from the weak convergence \eqref{wcx}, it is clear that
\begin{equation*}
    \int_0^T e^{2sA^\ast}\tau^7 (\vepx-\varphi_x)(t,1)v(t)\,dt =\int_0^Te^{s\hat{A}}\tau^{-1/2} (\vepx-\varphi_x)(t,1)w(t)v(t)\,dt\to 0,
\end{equation*}
which proves the result.
\end{proof}

\begin{lem}\label{lem5.2}
For each $0<\ep<1$, the linear operator $L_\ep:L^2(0,T;H^2_\alpha)\longrightarrow \R$, defined by 
\[L_\ep v:=-\intw \fep v\,dxdt=\frac{1}{\ep^3}\intw e^{2sA^\ast} \tau^7\vep v\, dxdt, \]
is bounded and weakly converges, up to a subsequence, as $\ep\to 0^+$, to the linear operator $L: L^2(0,T;H^2_\alpha)\longrightarrow \R$ defined by 
\[L v:=\frac{1}{3}\int_0^T e^{2sA^\ast} \tau^7\varphi_x(t,1) v_x(t,1)\,dt. \]
\end{lem}

\begin{proof}
Firstly, let us prove the boundedness of $L_\ep$. Given $v\in L^2(0,T;H^2_\alpha)$, from \eqref{3.3*}, we can see that
\begin{equation}\label{4.1}
    |L_\ep v|\leq \intw |e^{-sA^\ast} \tau^{-7/2}\fep| |e^{sA^\ast} \tau^{7/2}v|\,dxdt \leq  \frac{C}{\ep^{3/2}}\left(\intw e^{2sA^\ast}\tau^7 v^2\,dxdt\right)^{1/2}\n{u_0}{L^2(0,1)}.
\end{equation}

In order to estimate the integral in the last inequality, for $s\in [1-\ep,1]$, first note that
\[\left| \int_x^1\int_s^1 {v_{xx}(t,r)}\,drds\right|\leq \int_x^1\int_s^1 |{v_{xx}(t,r)}|\,drds\leq (1-x)\int_{1-\ep}^1|{v_{xx}(t,r)}|\,dr.\]

As a consequence of Jensen's inequality, we have
\[\left| \int_x^1\int_s^1{v_{xx}(t,r)}\,drds\right|^2 \leq (1-x)^2\left(\int_{1-\ep}^1|{v_{xx}(t,r)}|\,dr\right)^2 \leq \ep(1-x)^2\int_{1-\ep}^1|{v_{xx}(t,r)}|^2\,dr.\]

Now, let us set $a=1-\ep$ in \eqref{eq.v}. For any $x\in [1-\ep,1)$, this last inequality gives us 
\begin{align*}
|v(t,x)|^2 & \leq  \left(|1-x||v_x(t,1)|+\left|\int_x^1\int_s^1 {v_{xx}(t,r)}\,drds \right|\right)^2\\
    &  \leq 2  \left(|1-x|^2|v_x(t,1)|^2+\left|\int_x^1\int_s^1 {v_{xx}(t,r)}\,drds \right|^2\right)\\
    & \leq  2\left((1-x)^2|v_x(t,1)|^2+\ep(1-x)^2\int_{1-\ep}^1|{v_{xx}}(t,r)|^2dr\right).
\end{align*}

Hence, using this estimate in \eqref{4.1}, we have
\begin{align*}
 |L_\ep v|^2\  & \begin{multlined}[t]
  \leq \frac{C}{\ep^3}\left(\intw e^{2sA^\ast}\tau^7(1-x)^2 |v_x(t,1)|^2\,dxdt\right)\\
  \quad +\frac{C}{\ep^2}\left(\intw e^{2sA^\ast}\tau^7 (1-x)^2\left[\int_{1-\ep}^1|{v_{xx}}(t,r)|^2\,dr\right]dxdt\right)
\end{multlined}\\
 &=:I_1+I_2.  
    \end{align*}
Recalling that $e^{2sA^\ast}\tau^7$ is bounded, we can apply Corollary \ref{abalfa} to get
\begin{equation*}
    I_1\leq \frac{C}{\ep^3}\intw (1-x)^2 |v_x(t,1)|^2\,dxdt\leq C \int_0^T|v_x(t,1)|^2\,dt\leq C_\alpha\n{v}{L^2(0,T;H^2_\alpha)}^2.
\end{equation*}

Besides, in order to estimate $I_2$, we note that we can consider $\ep\leq 1/2$, since, in the end, it will go to zero. Thus, from \eqref{est.uxx}, we conclude that
\begin{align*}
    I_2 & \leq 4C \left(\frac{1}{(1-\ep)^{2\alpha}}+\left(\frac{\alpha}{1-\ep}\right)^2\right) \intw e^{2sA^\ast}\tau^7\n{v(t)}{H^2_\alpha}^2\,dxdt\\
    &  \leq C\ep \left(\frac{1}{(1-\ep)^{2\alpha}}+\left(\frac{\alpha}{1-\ep}\right)^2\right) \int_0^T e^{2sA^\ast}\tau^7\n{v(t)}{H^2_\alpha}^2\,dt\\
    & \leq C_\alpha \n{v}{L^2(0,T;H^2_\alpha)}^2.
\end{align*}
It means that $\{ L_{\varepsilon} ; 0<\varepsilon <1/2\}$ is an uniformly bounded family in  $\left(L^2(0,T;H^2_\alpha)\right)'$.
Therefore, there exists $L\in \left(L^2(0,T;H^2_\alpha)\right)'$ such that, up to a subsequence,
\[L_\ep \rightharpoonup L \text{ weakly in } L^2(0,T;H^2_\alpha) \text{ as } \ep\to 0^+. \]

It remains to prove that
\begin{equation}\label{L}
L v=\frac{1}{3}\intw e^{2sA^\ast} \tau^7\vepx(t,1) v_x(t,1)\, dxdt.    
\end{equation}
Indeed, if we define $V(t,x)=\frac{1}{1-x}\int_x^1\int_s^1 v_{xx}(t,r)\,drds$ we can rewrite  \eqref{eq.v} as
\begin{equation*}
    v(t,x)=-(1-x)v_x(t,1)+(1-x)V(t,x).
\end{equation*}
Thus, 
\begin{align*}
    L_\ep v & =-\frac{1}{\ep^3}\intw e^{2sA^\ast} \tau^7\vep(1-x)v_x(t,1)\, dxdt+\frac{1}{\ep^3}\intw e^{2sA^\ast} \tau^7\vep(1-x)V(t,x)\, dxdt\\
    & =:A_\ep+B_\ep.
    \end{align*}
In the following, we will prove that $B_\ep \to 0$ and $A_\ep$ converges to the right hand side of \eqref{L}. In fact, since $e^{sA^\ast}\tau^7$ is bounded, we can see that
\begin{align*}
    |B_\ep|& \leq \left(\frac{1}{\ep^3}\intw e^{2sA^\ast}\tau^7|\vep|^2\,dxdt \right)^{1/2}\left(\frac{1}{\ep^3}\intw e^{2sA^\ast}\tau^7(1-x)^2V(t,x)^2 \,dxdt\right)^{1/2}\\
    & \leq C\|\vep\|_P\left(\frac{1}{\ep^3}\intw (1-x)^2V(t,x)^2\,dxdt \right)^{1/2}\\
    & \leq C\|\vep\|_P\sup_{x\in [1-\ep,1]}\n{V(\cdot,x)}{L^2(0,T)}\to 0,\ \text{ as } \ep\to 0^+,
\end{align*}
as proved in \eqref{lim-sup}.

Now, let us obtain the convergence of $A_\ep$. Since $\vep(t,\cdot)\in H^2(1-\ep,1)$, we observe that
\begin{align*}
    A_\ep = & \frac{1}{\ep^3}\intw e^{2sA^\ast} \tau^7(1-x)v_x(t,1)\left(\int_x^1\vepx(t,s)ds\right)\, dxdt\\
    = &\begin{multlined}[t]
    \frac{1}{\ep^3}\intw e^{2sA^\ast} \tau^7(1-x)v_x(t,1)\left(\int_x^1\vepx(t,s)-\vepx(t,1)ds\right)\, dxdt \\
    +\frac{1}{\ep^3}\intw e^{2sA^\ast} \tau^7(1-x)^2v_x(t,1)\vepx(t,1)\,dxdt
    \end{multlined}\\
    =: & A_\ep^1+A_\ep^2.
    \end{align*}
Hence, from Proposition \ref{prop4.1},
\begin{equation*}
    A_\ep^2=\frac{1}{3}\int_0^Te^{2sA^\ast} \tau^7v_x(t,1)\vepx(t,1)\,dt\to \frac{1}{3}\int_0^Te^{2sA^\ast} \tau^7v_x(t,1)\vep(t,1)\,dt.
\end{equation*}
In order to prove $A_\ep^1\to 0$, we firstly recall the definition of $\bar{A}$ given in \eqref{functions2} and note that $e^{2s\bar{A}}\tau^{15}$ is bounded. In the following, since $\vepx\in H^1(1-\ep,1)$, we can apply \eqref{TFC}. Finally, we use Corollary \ref{abalfa},  inequality \eqref{est.uxx} and Proposition \ref{prop4.1} to obtain
\begin{align*}
    |A^1_\ep|\leq & \frac{1}{\ep^3}\intw e^{2sA^\ast} \tau^7(1-x)|v_x(t,1)|\left(\int_x^1\int_s^1|\vepxx(t,\xi)|\,d\xi ds\right)\, dxdt \\
     = &\begin{multlined}[t]
    \frac{1}{\ep^3}\left(\intw e^{2s\bar{A}}\tau^{15}|v_x(t,1)|^2\,dxdt\right)^{1/2}\\  \quad \quad\left(\intw e^{2s\hat{A}}\tau^{-1}(1-x)^2\left(\int_x^1\int_s^1|\vepxx(t,\xi)|\,d\xi ds\right)^2\,dxdt\right)^{1/2}
    \end{multlined}\\
    \leq  & \frac{1}{\ep^3}\left(\ep\int_0^T |v_x(t,1)|^2\,dt\right)^{1/2}\left(\intw e^{2s\hat{A}}\tau^{-1}(1-x)^3\left[\int_x^1(1-s)\int_{s}^1|\vepxx(t,\xi)|^2\,d\xi ds\right]\, dxdt\right)^{1/2}\\
    \leq &  \frac{C_\alpha\ep^{1/2}}{\ep^3}\n{v}{L^2(0,T;H^2_\alpha)}\left(\intw e^{2s\hat{A}}\tau^{-1}(1-x)^4\left[\int_{1-\ep}^1\int_{1-\ep}^1|\vepxx(t,\xi)|^2\,d\xi ds\right]\,dxdt\right)^{1/2}\\
     \leq  &   \frac{C_\alpha\ep^{1/2}}{\ep^3}\n{v}{L^2(0,T;H^2_\alpha)} \left(\frac{1}{(1-\ep)^{2\alpha}}+\frac{\alpha^2}{(1-\ep)^{2+\alpha}}\right)^{1/2}\left(\intw e^{2s\hat{A}}\tau^{-1}(1-x)^4\ep^2\n{\vep(t)}{H_\alpha^2}^2\, dxdt\right)^{1/2}\\
    \leq      & \frac{C_\alpha\ep^{1/2}}{\ep^3}\n{v}{L^2(0,T;H^2_\alpha)} \left(\frac{1}{(1-\ep)^{2\alpha}}+\frac{\alpha^2}{(1-\ep)^{2+\alpha}}\right)^{1/2}\left(\frac{\ep^7}{7}\intw e^{2s\hat{A}}\tau^{-1}\n{\vep(t)}{H_\alpha^2}^2\, dxdt\right)^{1/2}\\
\leq &  C_\alpha\ep\n{v}{L^2(0,T;H^2_\alpha)} \left(\frac{1}{(1-\ep)^{2\alpha}}+\frac{\alpha^2}{(1-\ep)^{2+\alpha}}\right)^{1/2} \left(\intw e^{2s\hat{A}}\tau^{-1}\n{\vep(t)}{H_\alpha^2}^2\, dxdt\right)^{1/2}\\
\leq &   C_\alpha\ep\n{v}{L^2(0,T;H^2_\alpha)} \left(\frac{1}{(1-\ep)^{2\alpha}}+\frac{\alpha^2}{(1-\ep)^{2+\alpha}}\right)^{1/2} \n{u_0}{L^2(0,1)}\to 0, \text{ as } \ep\to 0^+.
\end{align*}
This concludes the result.
\end{proof}


Finally, we are ready to obtain our main result (Theorem \ref{main}) as a consequence of the previous lemma.

\begin{proof}[Proof of Theorem \ref{main}]

Recall $(\yep,\fep)=(e^{2sA}\mathcal{L}^*\hat{\varphi}_\varepsilon,-\varepsilon^{-3}e^{2sA^*}\tau^7\hat{\varphi}_\varepsilon)$ is the solution of the internal null controllability problem \eqref{dcp} in the sense of Definition \ref{transp}. From \eqref{sol.dcp}, for every $(F,z^T)\in L^1(Q)\times L^2(0,1)$, we have
\begin{equation}\label{soltrans}
    \int_0^T\int_0^1 \yep F\,dx\,dt=-\frac{1}{\ep^3}\int_0^T\int_{1-\varepsilon}^1e^{2sA^*}\tau^7\hat{\varphi}_\varepsilon z\,dx\,dt+\int_0^1u_0(x)z(x,0)\,dx,
\end{equation}
where $z$ is the associated solution to \eqref{adj-transp}. The estimate  \eqref{3.3*} gives us $\yep  \rightharpoonup u$ in $L^2(Q)$, which can be combined with Lemma \ref{lem5.2}, in order to pass to the limit in \eqref{soltrans}, following
\begin{equation*}
    \int_0^T\int_0^1 uF\,dx\,dt=\frac{1}{3}\int_0^Te^{2sA^*}\tau^7\varphi_x(t,1) z_x(t,1)\,dt+\int_0^1u_0(x)z(x,0)\,dx.
\end{equation*}
In other words, the pair $(u,\frac{1}{3}e^{2sA^*}\tau^7\varphi_x(\cdot,1))$ solves the boundary null controllability problem \eqref{bcp}. \end{proof}

\section{Further extensions and open questions}\label{aberto}

In this section we will present some additional comments on the main result and the open questions left for future work. 

\subsection{A more general degenerate operator}

In most works published on null controllability of degenerate parabolic equations, the authors deal with the following base system
\textit{distributed control problem:}
	\begin{equation}
	\begin{cases}
	u_t- \displaystyle \left(a(x) u_x \right)_x=h\chi_{\omega_\ep}, & (t,x)\in Q,\\
	u(t,1)=0,& \text{ in } (0,T),\\
	\begin{cases*}
	u(t,0)=0, & \text{ (Weak) }\\
	\text{or}\\
	(a(x) u_x)(t,0)=0,&  \text{ (Strong) }
	\end{cases*},&  t\in (0,T), \\
	u(0,x)=u_0(x),
	\end{cases}
	\label{gdcp}
	\end{equation}
where $a$ is a positive continuous function on $[0,1]$ that satisfies the following hypothesis
\begin{eqnarray*}
(Weak): \  \left\{\begin{array}{ccc}(i) \ \ a\in C([0,1])\cap C^1((0,1]), \ a>0 \ \mbox{in} \ (0,1] \ \mbox{and} \ a(0)=0;\\ (ii) \ \ \exists K\in[0,1) \ \mbox{such that} \ xa'(x)\leq Ka(x) \ \forall x\in[0,1].\hspace{0.6cm}\end{array}\right.\hspace{1.5cm}\\ \\
(Strong): \ \left\{\begin{array}{ccc}(i) \ \ a\in C^1([0,1]), \ a>0 \ \mbox{in} \ (0,1] \ \mbox{and} \ a(0)=0;\hspace{3cm}\\ (ii) \ \ \exists K\in[1,2) \ \mbox{such that} \ xa'(x)\leq Ka(x) \ \forall x\in[0,1];\hspace{2cm}\\ (iii) \ \ \left\{\begin{array}{cc}\exists \kappa\in(1,K]; x\to\frac{a(x)}{x^\kappa} \ \mbox{is nondecreasing near 0, if} \ K>1;\\ \exists \kappa\in(0,1); x\to\frac{a(x)}{x^\kappa} \ \mbox{is nondecreasing near 0, if} \ K=1.\end{array}\right. \end{array}\right.
\end{eqnarray*}

In this work we deal with the particular case where $a(x)=x^\alpha$, but there are no significant changes to adapt our computations to deal with problem \eqref{gdcp}. So, Theorem \ref{main} can be extended naturally to the problem \eqref{gdcp}.

\subsection{The problem in two dimensional space}

There is a lot of ways to extend the system \eqref{dcp} in dimension 2. One of the simplest one is to consider
\begin{equation}
\label{prob1}
\left\{\begin{array}{lll}u_t-\text{\rm div\,}(A\nabla u)  =  h\chi_{\omega}\ & \text{in}&   Q,\\ 
B.C. &\mbox{on}& \Sigma,\\u(\cdot,0)= u_0  & \text{in} &  \Omega,\\ \end{array}\right.
\end{equation} 
where $\Omega:=(0,1)\times(0,1)$, $\Gamma:=\partial\Omega$, $T>0$, $Q:=\Omega\times(0,T)$, $\Sigma:=\Gamma\times(0,T)$, $\omega\subset\Omega$ is a non-empty open set , $\chi_\omega$ is the corresponding characteristic function, $h\in L^2(Q)$, $u_0\in L^2(\Omega)$, and $A:\overline{\Omega}\longrightarrow  M_{2\times 2}(\mathbb{R})$ is the matrix-valued function given by 
\[
A(x)=diag(x_1^{\alpha_1},x_2^{\alpha_2})
\]
In \eqref{prob1}, the boundary conditions are
$$B.C.:=\left\{\begin{array}{llll}u=0\ \ \mbox{on}\ \  \Sigma   &\mbox{if}& \alpha_1,\ \alpha_2\in[0,1),\\ u=0 \ \  \mbox{on} \ \ \Sigma_{3,4} \ \  \mbox{and} \ \ (A\nabla u)\nu=0\ \ \mbox{on} \ \ \Sigma_{1,2}&\mbox{if}& \alpha_1,\ \alpha_2\in[1,2),\\ u=0 \ \ \mbox{on} \ \ \Sigma_{1,3,4} \ \ \mbox{and} \ \ (A\nabla u)\nu=0 \ \ \mbox{on}  \ \ \Sigma_2& \mbox{if}& \alpha_1\in[0,1) \ \  \mbox{and} \ \ \alpha_2\in[1,2),\\ u=0 \ \ \mbox{on} \ \ \Sigma_{2,3,4} \ \ \mbox{and} \ \ (A\nabla u)\nu=0 \ \ \mbox{on}  \ \ \Sigma_1& \mbox{if}& \alpha_1\in[1,2) \ \  \mbox{and} \ \ \alpha_2\in[0,1),\\  \end{array}\right.$$
with $\alpha=(\alpha_1,\alpha_2)\in[0,2)\times[0,2)$, $\Sigma_{i,j,l}:=(\Gamma_i\cup\Gamma_j\cup\Gamma_l)\times(0,T)$, $\nu=\nu(x)$ being the outward unit normal to $\Omega$ at the point $x \in \Gamma$ and, finally, \begin{equation*}\label{boundary}\Gamma_1:=\{0\}\times[0,1], \ \Gamma_2:=[0,1]\times\{0\}, \ \Gamma_3:=\{1\}\times[0,1] \text{and} \ \Gamma_4:=[0,1]\times\{1\}.
\end{equation*}

In \cite{araruna2019carleman} the authors established the null controllability of \eqref{prob1} by using of a internal Carleman estimate. To adapt the ideas in this work to deal with system \eqref{prob1}, we need a boundary Carleman estimate and, to our best knowledge, a such estimate does not exist. However, this does not mean that it is  impossible question to solve. 

\subsection{The boundary control acting at 0}

A natural question that this work left is whether the Theorem \ref{main} works with the boundary control $g$ acting at $x=0$ instead of $x=1$. This is a non trivial question. Firstly, the boundary null controllability problem \eqref{bcp} with the control $g$ acting at $x=0$ instead of $x=1$ seems to be much more complicated to deal with. Indeed, to our best knowledge, \eqref{bcp} with the control $g$ acting at $x=0$, was only solved for the weak degenerate case in \cite{gueye2014exact} by using of the momentum method instead of Carleman estimates. One of the main difficulties to deal, in the strong degenerate case, is the fact that there are no trace results at $x=0$ when $\alpha\geq1$. So, to be hopeful about using the ideas of this work, on order to get a similar result with the control acting at $x=0$ instead of $x=1$, the best bet would  be only the consideration of the weak degenerate case. However, even in this situation, to use the ideas of this work, it is expected a boundary Carleman estimate with the control acting $x=0$, which may be a hard task.

\appendix
\renewcommand{\thesection}{\Alph{section}} 

\section{Appendix}\label{app}
At this place, we will complete the proof of Proposition \ref{car-frnt}. To be more precise, we will give detailed explanations about \eqref{uu}. As in Section \ref{sec-pre}, let us consider
\[w(t,x):=e^{s\varphi(t,x)}v(t,x),\]
for each $(0,T)\times (0,1)$, where $v$ is a solution of \eqref{adj-jlr}. Putting 
\[Lv:=v_t+(x^{\alpha} v_x)_x \ \ \mbox{ and }\ \  L_sw:=e^{s\varphi}L(e^{-s\varphi}w), \emph{ with } s>0,
\]
we notice that $w$ satisfies
\begin{equation}\lnum \label{prob3.4.1}
 \left\{\begin{array}{ll}
L_sw=e^{s\varphi} F(t,x), & (t,x)\in (0,T)\times (0,1), \\
w(t,0)=w(t,1)=0, & t\in (0,T),\\
w(T,x)=w(0,x)=0, & x\in (0,1).
\end{array}\right.   \end{equation}
Besides, decomposing $L_sw$ as
$$L_sw:=L_s^+w+L_s^-w,$$
where
\[
L_s^+w:=(x^{\alpha} w_x)_x-s\varphi_tw+s^2 x^{\alpha} \varphi_x^2 w
\]
and
\[
L_s^-w:=w_t-2sx^{\alpha} \varphi_x w_x-s(x^{\alpha} \varphi_x)_x w,
\]
we certainly get
$$\|L_s^+w\|^2+\|L_s^-w\|^2+2\langle L_s^+w,L_s^-w \rangle=\|F e^{s\varphi}\|^2.$$
Under these notations, we are ready to obtain the next result.

\begin{lem}\label{prop3.3}
Let $T>0$ be given. Then, there exist $C>0$ and $s_0>0$, both independent of $w$, such that
\begin{multline}\label{uu+}
 \intq(s\theta)^{-1}(w_t^2+|(x^{\alpha} w_x)_x|^2)\,dxdt   +\intq \left(s\theta x^{\alpha} w_x^2+ (s\theta)^{3} x^{2-\alpha} w^2\right)\,dxdt\\
 \leq C\left(\intq e^{2s\varphi}|F|^2+s \int_0^T \theta(t)w_x^2(t,1)\,dxdt    \right).
\end{multline}
for any $s\geq s_0$. In particular, \eqref{uu} holds. 
\end{lem}

\begin{proof} 

In \cite{alabau2006carleman}, 
the authors have obtained the following  estimate involving the scalar product $\langle L_s^+w,L_s^-w \rangle$:
\begin{equation}\label{myeq21}
\intq \left(s\theta x^{\alpha} w_x^2+(s\theta)^3 x^{2-\alpha} w^2\right)dxdt
\leq C\bigg(2\langle L_s^+w,L_s^-w \rangle + \left.s\int_0^T \theta(t)w_x^2(t,1)dt\    \right).
\end{equation}
In order to obtain the desired inequality \eqref{uu+}, we will estimate
\begin{equation}\label{myeq16} 
\intq(s\theta)^{-1}|(x^{\alpha} w_x)_x|^2 dxdt  \mbox{ and } \intq(s\theta)^{-1}|w_t|^2 dxdt  
\end{equation}
in terms of $\|L_s^+w\|^2$ and $\|L_s^-w\|^2$.
Since
\begin{align*}
|(x^{\alpha} w_x)_x|^2\leq 3(|L_s^+w|^2+s^2|\varphi_t|^2|w|^2+s^4x^{2\alpha}\varphi_x^4w^2)
\end{align*}
we get
\begin{align}\label{myeq10}
  \intq(s\theta)^{-1}|(x^{\alpha} w_x)_x|^2 dxdt  
  & \leq 3\intq(s\theta)^{-1}|L_s^+ w|^2 dxdt \nonumber\\
 & +3\intq(s\theta)^{-1}s^2|\varphi_t|^2|w|^2dxdt \nonumber\\
&  +3\intq(s\theta)^{-1}s^4 x^{2\alpha}\varphi_x^4w^2 dxdt
\nonumber\\
&=:3\intq(s\theta)^{-1}|L_s^+ w|^2 dxdt + \mathcal I + \mathcal J.
\end{align}
Let us estimate $\mathcal I$ and $\mathcal J$ separately. Firstly, since $|\theta^{-1}\theta_t^2\psi^2|\leq C\theta^{2}$ and $s\leq Cs^2$ for some $C>0$, Hardy-Poincar\'{e} inequality yields

\begin{align}  \label{myeq12}
\mathcal I
&=3\intq s|\theta^{-1}\theta_t^2\psi^2|w^2 dxdt \nonumber \\
&\leq C\intq (s\theta)^2w^2 dxdt\nonumber \\
& =C\intq \left((s\theta)^{1/2} x^{\frac{\alpha -2}{2}} w\right)\left((s\theta)^{3/2} x^{\frac{2-\alpha}{2}} w\right)dxdt \nonumber \\
& \leq C\left(\intq (s\theta) x^{\alpha -2} w^2 dxdt
+\intq(s\theta)^{3} x^{2-\alpha} w^2\right) dxdt \nonumber \\
&\leq C\left(\intq (s\theta) x^{\alpha} w_x^2 dxdt +\intq(s\theta)^{3} x^{2-\alpha} w^2 dxdt \right).
\end{align}
Secondly, to deal with $\mathcal J$, we recall $\varphi_x =\theta x^{1-\alpha}$ and $x^{2-\alpha} \leq 1$. Thus,

\begin{equation}\label{myeq14}
\mathcal J
=\intq (s\theta)^{-1} s^4 x^{2\alpha} (\theta ^4 x^{4-4\alpha}) w^2 dxdt
\leq C  \intq(s\theta)^{3} x^{2-\alpha} w^2 dxdt.
\end{equation}
As a  result, we can combine (\ref{myeq10}), (\ref{myeq12}) and (\ref{myeq14}) to obtain
\begin{align}\label{myeq15}
 & \intq(s\theta)^{-1}|(x^{\alpha} w_x)_x|^2 dxdt   \nonumber\\
 & \leq C\left(\|L_s^+w\|^2
 +\intq (s\theta ) x^{\alpha} w_x^2 dxdt  + \intq (s\theta)^3 x^{2-\alpha} w^2 dxdt  \right).
\end{align}

Next, we will focus on the second integral mentioned in \eqref{myeq16}.
Notice that 
$$w_t^2\leq 3\left(|L_s^-w|^2+4s^2 x^{2\alpha} \varphi_x^2w_x^2+s^2 x^{2\alpha} \varphi^4w^2\right)$$
implies
\begin{align}\label{myeq17}
 \intq(s\theta)^{-1}w_t^2 dxdt 
 &= 3\intq (s\theta)^{-1}|L_s^-w|^2 dxdt \nonumber \\
& +12\intq (s\theta)^{-1}s^2 x^{2\alpha} \varphi_x^2 w_x^2 dxdt \nonumber\\
& +3\intq (s\theta)^{-1} s^2 x^{2\alpha} \varphi^4 w^2 dxdt \nonumber\\
&=:3\intq (s\theta)^{-1}|L_s^-w|^2 dxdt
+\mathcal K
+\mathcal M .
\end{align}
Clearly, using twice the relation $x^2 \leq x^{\alpha}$, we have
\begin{equation}\label{myeq18}
\mathcal K 
=12\intq(s\theta)^{-1}s^2 x^{2\alpha} (\theta ^2 x^{2-2\alpha}) w_x^2 dxdt
\leq C\intq s\theta x^{\alpha} w_x^2 dxdt  
\end{equation}
and
\begin{equation}\label{myeq19}
\mathcal M =\intq (s\theta)^{-1}s^2 x^{2\alpha} (\theta ^4 x^{4-4\alpha }) w^2 dxdt  
\leq C\intq (s\theta)^{3} x^{2-\alpha} w^2 dxdt,  
\end{equation}
following
\begin{align}\label{myeq20}
& \intq(s\theta)^{-1}w_t^2 dxdt  \nonumber\\
&\leq C\left(\|L_s^-w\|^2+\intq s\theta x^{\alpha} w_x^2 dxdt  + \intq (s\theta)^{3} x^{2-\alpha} w^2 dxdt\right).
\end{align}
As in \eqref{myeq15}, we have also used in \eqref{myeq20} that $(s\theta )^{-1}$ is bounded.

Therefore, from \eqref{myeq21}, \eqref{myeq15} and \eqref{myeq20}, we conclude that
\begin{align*}
& \intq(s\theta)^{-1}(w_t^2+|(x^{\alpha} w_x)_x|^2)dxdt   
+\intq (s\theta x^{\alpha} w_x^2+ (s\theta)^{3} x^{2-\alpha} w^2 )dxdt \\
&\leq C\left(\|L_s^+w\|^2+\|L_s^-w\|^2+\intq (s\theta x^{\alpha} w_x^2   
+ (s\theta)^{3} x^{2-\alpha} w^2)dxdt  \right) \\
& \leq C\left(\|L_s^+w\|^2+\|L_s^-w\|^2+2\langle L_s^+w,L_s^-w \rangle + s\int_0^T \theta(t)w_x^2(t,1) dt \right) \\
& \leq C\left(\intq e^{2s\varphi}|F|^2 dxdt
+s\int_0^T \theta(t)w_x^2(t,1)dt \right),
\end{align*}
as expected.
\end{proof}

\bibliography{references}

\begin{thebibliography}{24}
\expandafter\ifx\csname natexlab\endcsname\relax\def\natexlab#1{#1}\fi
\providecommand{\url}[1]{\texttt{#1}}
\providecommand{\href}[2]{#2}
\providecommand{\path}[1]{#1}
\providecommand{\DOIprefix}{doi:}
\providecommand{\ArXivprefix}{arXiv:}
\providecommand{\URLprefix}{URL: }
\providecommand{\Pubmedprefix}{pmid:}
\providecommand{\doi}[1]{\href{http://dx.doi.org/#1}{\path{#1}}}
\providecommand{\Pubmed}[1]{\href{pmid:#1}{\path{#1}}}
\providecommand{\bibinfo}[2]{#2}
\ifx\xfnm\relax \def\xfnm[#1]{\unskip,\space#1}\fi
\bibitem[{Alabau-Boussouira et~al.(2006)Alabau-Boussouira, Cannarsa \&
  Fragnelli}]{alabau2006carleman}
\bibinfo{author}{Alabau-Boussouira, F.}, \bibinfo{author}{Cannarsa, P.}, \&
  \bibinfo{author}{Fragnelli, G.} (\bibinfo{year}{2006}).
\newblock \bibinfo{title}{Carleman estimates for degenerate parabolic operators
  with applications to null controllability}.
\newblock {\it \bibinfo{journal}{Journal of Evolution Equations}\/},  {\it
  \bibinfo{volume}{6}\/}, \bibinfo{pages}{161--204}.
\bibitem[{Araruna et~al.(2019)Araruna, Ara{\'u}jo \&
  Fern\'andez-Cara}]{araruna2019carleman}
\bibinfo{author}{Araruna, F.~D.}, \bibinfo{author}{Ara{\'u}jo, B. S.~V.}, \&
  \bibinfo{author}{Fern\'andez-Cara, E.} (\bibinfo{year}{2019}).
\newblock \bibinfo{title}{Carleman estimates for some two-dimensional
  degenerate parabolic pdes and applications}.
\newblock {\it \bibinfo{journal}{SIAM Journal on Control and Optimization}\/},
  {\it \bibinfo{volume}{57}\/}, \bibinfo{pages}{3985--4010}.
\bibitem[{Araruna et~al.({2018})Araruna, Araújo \&
  Fernandez-Cara}]{araruna2018stackelberg}
\bibinfo{author}{Araruna, F.~D.}, \bibinfo{author}{Araújo, B. S.~V.}, \&
  \bibinfo{author}{Fernandez-Cara, E.} (\bibinfo{year}{{2018}}).
\newblock \bibinfo{title}{{Stackelberg-Nash null controllability for some
  linear and semilinear degenerate parabolic equations}}.
\newblock {\it \bibinfo{journal}{Mathematics of Control Signals and
  Systems}\/},  {\it \bibinfo{volume}{{30}}\/}.
\bibitem[{Boutaayamou et~al.({2018})Boutaayamou, Fragnelli \&
  Maniar}]{Boutaayamou2018carleman}
\bibinfo{author}{Boutaayamou, I.}, \bibinfo{author}{Fragnelli, G.}, \&
  \bibinfo{author}{Maniar, L.} (\bibinfo{year}{{2018}}).
\newblock \bibinfo{title}{{Carleman estimates for parabolic equations with
  interior degeneracy and Neumann boundary conditions}}.
\newblock {\it \bibinfo{journal}{{Journal d'Analyse Math\'ematique}}\/},  {\it
  \bibinfo{volume}{{135}}\/}, \bibinfo{pages}{{1--35}}.
\bibitem[{Brezis(1983)}]{Brezis}
\bibinfo{author}{Brezis, H.} (\bibinfo{year}{1983}).
\newblock {\it \bibinfo{title}{Analyse Fonctionnelle: Théorie et
  Applications}\/}.
\newblock \bibinfo{address}{Paris}: \bibinfo{publisher}{Masson}.
\bibitem[{Campiti et~al.(1998)Campiti, Metafune \&
  Pallara}]{campiti1998degenerate}
\bibinfo{author}{Campiti, M.}, \bibinfo{author}{Metafune, G.}, \&
  \bibinfo{author}{Pallara, D.} (\bibinfo{year}{1998}).
\newblock \bibinfo{title}{Degenerate self-adjoint evolution equations on the
  unit interval}.
\newblock In {\it \bibinfo{booktitle}{Semigroup Forum}\/} (pp.
  \bibinfo{pages}{1--36}).
\newblock \bibinfo{organization}{Springer} volume~\bibinfo{volume}{57}.
\bibitem[{Cannarsa \& De~Teresa(2009)}]{cannarsa2009controllability}
\bibinfo{author}{Cannarsa, P.}, \& \bibinfo{author}{De~Teresa, L.}
  (\bibinfo{year}{2009}).
\newblock \bibinfo{title}{Controllability of 1-d coupled degenerate parabolic
  equations.}
\newblock {\it \bibinfo{journal}{Electronic Journal of Differential
  Equations}\/},  {\it \bibinfo{volume}{2009}\/}, \bibinfo{pages}{1--21}.
\bibitem[{Cannarsa \& Fragnelli(2006)}]{cannarsa2006null}
\bibinfo{author}{Cannarsa, P.}, \& \bibinfo{author}{Fragnelli, G.}
  (\bibinfo{year}{2006}).
\newblock \bibinfo{title}{Null controllability of semilinear degenerate
  parabolic equations in bounded domains.}
\newblock {\it \bibinfo{journal}{Electronic Journal of Differential
  Equations}\/},  (pp. \bibinfo{pages}{1--20}).
\bibitem[{Cannarsa et~al.(2007)Cannarsa, Fragnelli \&
  Rocchetti}]{cannarsa2007null}
\bibinfo{author}{Cannarsa, P.}, \bibinfo{author}{Fragnelli, G.}, \&
  \bibinfo{author}{Rocchetti, D.} (\bibinfo{year}{2007}).
\newblock \bibinfo{title}{Null controllability of degenerate parabolic
  operators with drift}.
\newblock {\it \bibinfo{journal}{Networks \& Heterogeneous Media}\/},  {\it
  \bibinfo{volume}{2}\/}, \bibinfo{pages}{695}.
\bibitem[{Cannarsa et~al.(2008{\natexlab{a}})Cannarsa, Fragnelli \&
  Rocchetti}]{cannarsa2008controllability}
\bibinfo{author}{Cannarsa, P.}, \bibinfo{author}{Fragnelli, G.}, \&
  \bibinfo{author}{Rocchetti, D.} (\bibinfo{year}{2008}{\natexlab{a}}).
\newblock \bibinfo{title}{Controllability results for a class of
  one-dimensional degenerate parabolic problems in nondivergence form}.
\newblock {\it \bibinfo{journal}{Journal of Evolution Equations}\/},  {\it
  \bibinfo{volume}{8}\/}, \bibinfo{pages}{583--616}.
\bibitem[{Cannarsa et~al.(2002)Cannarsa, Martinez \&
  Vancostenoble}]{cannarsa2002nulle}
\bibinfo{author}{Cannarsa, P.}, \bibinfo{author}{Martinez, P.}, \&
  \bibinfo{author}{Vancostenoble, J.} (\bibinfo{year}{2002}).
\newblock \bibinfo{title}{Nulle contr{\^o}labilit{\'e} r{\'e}gionale pour des
  {\'e}quations de la chaleur d{\'e}g{\'e}n{\'e}r{\'e}es}.
\newblock {\it \bibinfo{journal}{Comptes rendus-M{\'e}canique}\/},  {\it
  \bibinfo{volume}{6}\/}, \bibinfo{pages}{397--401}.
\bibitem[{Cannarsa et~al.(2004)Cannarsa, Martinez \&
  Vancostenoble}]{cannarsa2004persistent}
\bibinfo{author}{Cannarsa, P.}, \bibinfo{author}{Martinez, P.}, \&
  \bibinfo{author}{Vancostenoble, J.} (\bibinfo{year}{2004}).
\newblock \bibinfo{title}{Persistent regional null contrillability for a class
  of degenerate parabolic equations}.
\newblock {\it \bibinfo{journal}{Communications on Pure \& Applied
  Analysis}\/},  {\it \bibinfo{volume}{3}\/}, \bibinfo{pages}{607}.
\bibitem[{Cannarsa et~al.(2008{\natexlab{b}})Cannarsa, Martinez \&
  Vancostenoble}]{cannarsa2008carleman}
\bibinfo{author}{Cannarsa, P.}, \bibinfo{author}{Martinez, P.}, \&
  \bibinfo{author}{Vancostenoble, J.} (\bibinfo{year}{2008}{\natexlab{b}}).
\newblock \bibinfo{title}{Carleman estimates for a class of degenerate
  parabolic operators}.
\newblock {\it \bibinfo{journal}{SIAM Journal on Control and Optimization}\/},
  {\it \bibinfo{volume}{47}\/}, \bibinfo{pages}{1--19}.
\bibitem[{Cannarsa et~al.(2016)Cannarsa, Martinez \&
  Vancostenoble}]{cannarsa2016global}
\bibinfo{author}{Cannarsa, P.}, \bibinfo{author}{Martinez, P.}, \&
  \bibinfo{author}{Vancostenoble, J.} (\bibinfo{year}{2016}).
\newblock {\it \bibinfo{title}{Global Carleman estimates for degenerate
  parabolic operators with applications}\/} volume \bibinfo{volume}{239, n.
  1133}.
\newblock \bibinfo{publisher}{American Mathematical Society}.
\bibitem[{Cannarsa et~al.(2005)Cannarsa, Martinez, Vancostenoble
  et~al.}]{cannarsa2005null}
\bibinfo{author}{Cannarsa, P.}, \bibinfo{author}{Martinez, P.},
  \bibinfo{author}{Vancostenoble, J.} et~al. (\bibinfo{year}{2005}).
\newblock \bibinfo{title}{Null controllability of degenerate heat equations}.
\newblock {\it \bibinfo{journal}{Advances in Differential Equations}\/},  {\it
  \bibinfo{volume}{10}\/}, \bibinfo{pages}{153--190}.
\bibitem[{Chaves-Silva et~al.(2020)Chaves-Silva, Puel \&
  Santos}]{chaves2020boundary}
\bibinfo{author}{Chaves-Silva, F.}, \bibinfo{author}{Puel, J.-P.}, \&
  \bibinfo{author}{Santos, M.} (\bibinfo{year}{2020}).
\newblock \bibinfo{title}{Boundary null controllability as the limit of
  internal controllability: The heat case.}
\newblock {\it \bibinfo{journal}{ESAIM: Control, Optimisation \& Calculus of
  Variations}\/},  {\it \bibinfo{volume}{26}\/}.
\bibitem[{Du({2019})}]{runmei2019null}
\bibinfo{author}{Du, R.} (\bibinfo{year}{{2019}}).
\newblock \bibinfo{title}{{Null controllability for a class of degenerate
  parabolic equations with the gradient terms}}.
\newblock {\it \bibinfo{journal}{{Journal of Evolution Equations}}\/},  {\it
  \bibinfo{volume}{{19}}\/}, \bibinfo{pages}{{585--613}}.
\bibitem[{El~Mustapha et~al.({2019})El~Mustapha, Mohamed \&
  Lahcen}]{mustapha2019algebraic}
\bibinfo{author}{El~Mustapha, A. B.~H.}, \bibinfo{author}{Mohamed, F.}, \&
  \bibinfo{author}{Lahcen, M.} (\bibinfo{year}{{2019}}).
\newblock \bibinfo{title}{{On Algebraic Condition for Null Controllability of
  Some Coupled Degenerate Systems}}.
\newblock {\it \bibinfo{journal}{{Mathematical Control and Related Fields}}\/},
   {\it \bibinfo{volume}{{9}}\/}, \bibinfo{pages}{{77--95}}.
\bibitem[{Fabre(1992)}]{fabre1992exact}
\bibinfo{author}{Fabre, C.} (\bibinfo{year}{1992}).
\newblock \bibinfo{title}{Exact boundary controllability of the wave equation
  as the limit of internal controllability}.
\newblock {\it \bibinfo{journal}{SIAM journal on control and optimization}\/},
  {\it \bibinfo{volume}{30}\/}, \bibinfo{pages}{1066--1086}.
\bibitem[{Fragnelli({2018})}]{fragnelli2018carleman}
\bibinfo{author}{Fragnelli, G.} (\bibinfo{year}{{2018}}).
\newblock \bibinfo{title}{{Carleman estimates and null controllability for a
  degenerate population model}}.
\newblock {\it \bibinfo{journal}{{Journal De Mathematiques Pures Et
  Appliquees}}\/},  {\it \bibinfo{volume}{{115}}\/},
  \bibinfo{pages}{{74--126}}.
\bibitem[{Gueye(2014)}]{gueye2014exact}
\bibinfo{author}{Gueye, M.} (\bibinfo{year}{2014}).
\newblock \bibinfo{title}{Exact boundary controllability of 1-d parabolic and
  hyperbolic degenerate equations}.
\newblock {\it \bibinfo{journal}{SIAM Journal on Control and Optimization}\/},
  {\it \bibinfo{volume}{52}\/}, \bibinfo{pages}{2037--2054}.
\bibitem[{Martinez \& Vancostenoble(2006)}]{martinez2006carleman}
\bibinfo{author}{Martinez, P.}, \& \bibinfo{author}{Vancostenoble, J.}
  (\bibinfo{year}{2006}).
\newblock \bibinfo{title}{Carleman estimates for one-dimensional degenerate
  heat equations}.
\newblock {\it \bibinfo{journal}{Journal Of Evolution Equations}\/},  {\it
  \bibinfo{volume}{6}\/}, \bibinfo{pages}{325--362}.
\bibitem[{Wang et~al.({2018})Wang, Zhou, Du \& Liu}]{wang2018carleman}
\bibinfo{author}{Wang, C.}, \bibinfo{author}{Zhou, Y.}, \bibinfo{author}{Du,
  R.}, \& \bibinfo{author}{Liu, Q.} (\bibinfo{year}{{2018}}).
\newblock \bibinfo{title}{Carleman estimate for solutions to a degenerate
  convection-diffusion equation}.
\newblock {\it \bibinfo{journal}{Discrete and Continuous Dynamical
  Systems-Series B}\/},  {\it \bibinfo{volume}{{23}}\/},
  \bibinfo{pages}{{4207--4222}}.
\bibitem[{Zuazua(1988)}]{zuazua1988controlabilite}
\bibinfo{author}{Zuazua, E.} (\bibinfo{year}{1988}).
\newblock \bibinfo{title}{Controlabilite exact interne de l'equation des
  ondes}.
\newblock {\it \bibinfo{journal}{Collection RMA}\/}, .

\end{thebibliography}
	
\end{document}